\documentclass[12pt]{amsart}
\usepackage{fullpage}
\usepackage{color}
\usepackage{pstricks,pst-node,pst-plot} %This loads the pstricks package
\usepackage{graphicx,psfrag} %,epstopdf}
\usepackage{color} %This loads the color package
\usepackage{tikz}
\usetikzlibrary{decorations.markings, hobby, knots}
\usepackage{pgffor}
\usepackage[pagebackref=true]{hyperref}
\usepackage[textsize=tiny]{todonotes}
\usepackage{subfigure}
\usepackage{enumitem}
\usepackage{mathtools}
\usepackage{verbatim}
\usepackage{amssymb}
\usepackage{amsmath}
\usepackage{amsthm}
\usepackage{float}
\usepackage{pinlabel}
\usepackage[dvipsnames]{xcolor}
\usepackage[normalem]{ulem}
\usepackage{multicol}
\usepackage{appendix}
\usepackage{ upgreek }

\usepackage{bm}
\usepackage{multirow}

\usepackage{perpage}	 %the perpage package
	
\definecolor{my_purple}{cmyk}{.5,.75,0,0}
\definecolor{my_darkblue}{cmyk}{1,1,0,0}
\definecolor{my_green}{cmyk}{1,0,1,0}
\definecolor{my_gray}{cmyk}{0,0,0,.5}
\definecolor{my_red}{cmyk}{0,1,1,0}
\definecolor{my_orange}{cmyk}{0,.35,.85,0}
\definecolor{my_cyan}{cmyk}{1,0,0,0}

\newtheorem{theorem}{Theorem}[section]
\newtheorem{corollary}[theorem]{Corollary}
\newtheorem{proposition}[theorem]{Proposition}
\newtheorem*{T:essentialconnected}{Theorem \ref{T:essentialconnected}}
\newtheorem*{C:essentialconnected}{Corollary \ref{C:essentialconnected}}

\theoremstyle{definition}
\newtheorem{definition}[theorem]{Definition}
\newtheorem{example}[theorem]{Example}

\theoremstyle{remark}
\newtheorem{remark}[theorem]{Remark}

\numberwithin{equation}{section}

\title{Moves connecting all essential continued fractions\\ of a rational number
}
\date{November 2026}

\email{cohenm@newpaltz.edu, tkindred@smith.edu, adlowrance@vassar.edu, pshanahan@lmu.edu, cvancott@usfca.edu}
\author{Moshe Cohen}
\address{State University of New York at New Paltz}
\author{Thomas Kindred}
\address{Smith College}
\author{Adam M. Lowrance}
\address{Vassar College}
\author{Patrick D. Shanahan}
\address{Loyola Marymount University}
\author{Cornelia A. Van Cott}
\address{University of San Francisco}

\begin{document}

\begin{abstract}
 Every non-integer rational number can be written in a variety of ways as a continued fraction with integer entries. Hirasawa and Teragaito described moves that reduce the length of a subtractive continued fraction without altering its value. They then showed that any two shortest expansions are related to each other by a finite repetition of an additional move. The latter result has geometric relevance in that it describes relationships among all minimal complexity spanning surfaces for a given 2-bridge knot or link.  We extend their result by expanding the class of continued fractions under consideration to those that correspond to $\pi_1$-essential spanning surfaces and introducing a new class of moves that, we prove, connects all such expansions of any rational number.
\end{abstract}

\maketitle

\section{Introduction}\label{S:intro1}

Every non-integer rational number can be written in a variety of ways as a continued fraction with integer entries. There are many natural questions one can ask about these continued fraction expansions, such as how short the shortest ones are, and how these expansions are related to each other.

These questions gain salience in light of a beautiful correspondence, first described by Schubert, between certain equivalence classes of rational numbers $\frac{p}{q}$ and 2-bridge knots and links $K_{p/q}$ \cite{Schubert_1954,Schubert_1956}. Further, work of Hatcher and Thurston tells us that under this correspondence, 
incompressible spanning surfaces for $K_{p/q}$ correspond to continued fraction expansions of $p/q$ in such a way that the complexity $\beta_1(F_\mathbf{b})$ of a surface equals the length $k$ of the corresponding continued fraction expansion $\mathbf{b}=r+[b_1,\dots,b_k]$ \cite{Hatcher-Thurston:1985}.  Moreover, $F_{\mathbf{b}}$ is orientable if and only if each term  $b_i$ of $\mathbf{b}$ is even, and $F_{\mathbf{b}}$ is $\pi_1$-essential if and only if each term is at least 2 in absolute value \cite[Theorem 1 (c)]{Hatcher-Thurston:1985}; in that case, we also call $\mathbf{b}$ {\it essential}. 

These correspondences enabled Hirasawa and Teragaito to determine the unoriented genus and crosscap number of $K_{p/q}$ by proving two theorems about continued fraction expansions \cite{Hirasawa_Teragaito:2006}.  First, they describe three moves---we call them 0-, $\pm1$-, and 232-moves in   \textsection\ref{subsec:HT}
---that take a continued fraction expansion of $\frac{p}{q}$ and transform it into a shorter one, and they prove that any sequence of these moves will terminate with a continued fraction of minimal length.  Second, they describe an additional move---which we call a $\pm2$-move---which takes a continued fraction expansion of $\frac{p}{q}$ and transforms it into another of the same length, and they prove that all minimal length expansions of $\frac{p}{q}$ are related by a sequence of these moves.  

While Hirasawa-Teragaito's results show that all expansions of $\frac{p}{q}$ are related by these four moves and their inverses, many of the expansions involved will be {\it inessential} in the sense that they have some term of $-1$, $0$, or $1$ and thus describe $\pi_1$-inessential surfaces.  In fact, Hirasawa-Teragaito's moves are not sufficient to relate all essential continued fraction expansions of $\frac{p}{q}$ directly---see Example \ref{Ex:1/3}.

In Definition \ref{def:n_move}, we generalize Hirasawa-Teragaito's $\pm2$-moves to a class of ``$\pm n$-moves'' for $2\leq n\in\mathbb{Z}$.
Whereas Hirasawa-Teragaito's $\pm2$-moves are inverses of one another this is not the case for our $\pm n$-moves: they have inverses which are a distinct class of moves. We call a $\pm n$-move $\mathbf{b}\to\mathbf{b}'$, or the inverse of such a move, {\it essential} if both $\mathbf{b}$ and $\mathbf{b}'$ are essential.  Our main result establishes that all essential expansions of $\frac{p}{q}$ are related by a sequence of essential $\pm n$ moves and their inverses.  We prove this by showing first that any essential expansion is related, in a particularly nice way, to the unique continued fraction expansion $\mathbf{b}_+$ of $\frac{p}{q}$ in which each term $b_i \ge 2$.

 \begin{theorem}
\label{T:essentialconnected}
Let $\frac{p}{q}\in\mathbb{Q}\setminus\mathbb{Z}$, and let $\mathbf{b}_+=r+[b_1,\dots,b_k]$ be the unique continued fraction expansion of $\frac{p}{q}$ in which each term $b_i\geq 2$. Every essential continued fraction expansion of $\frac{p}{q}$ can be obtained from $\mathbf{b}_+$ via a finite sequence of essential $n$-moves and essential inverse $-n$-moves for various $n\geq 2$.
    \end{theorem}

\begin{corollary}
\label{C:essentialconnected}
        Any two essential continued fraction expansions of $\frac{p}{q}\in \mathbb{Q}\setminus \mathbb{Z}$ are related by a sequence of essential $\pm n$-moves and their inverses for various $n\geq 2$. 
\end{corollary}

In our preceeding paper \cite{CohKinLowShaVan}, the authors obtained a simple formula for the length of the shortest continued fraction expansion of $\frac{p}{q}$ in terms of the continued fraction $\mathbf{b}_+$.  One can understand the main results of this paper as extensions on the same theme: the continued fraction $\mathbf{b}_+$  completely determines the structure of all essential continued fraction expansions of $\frac{p}{q}$.

In a follow-up \cite{CohKinLowShaVan:3} to the present paper, we will define {\it spanning surgery moves} of type $\pm n$, where $n$ is a positive integer, that correspond to $\pm n$-moves on continued fraction expansions.   Using this correspondence, we will establish a geometric version of Theorem~\ref{T:essentialconnected}, that any two $\pi_1$-essential spanning surfaces for $K_{p/q}$ are related by a finite sequence of surgery moves where each intermediate surface is also essential.  

For another perspective on Theorem \ref{T:essentialconnected}, let $\frac{p}{q}$ be a non-integer rational number, and consider the edge-weighted directed graph $G_{p/q}$ whose vertices correspond to the essential continued fraction expansions of $\frac{p}{q}$ and whose edges correspond to the essential $\pm n$-moves between these expansions.  Each edge is labeled accordingly with the value $\pm n$ and is directed as explained in Definition \ref{D:G_p/q}. Theorem~\ref{T:essentialconnected} asserts that the underlying undirected graph of $G_{p/q}$ is connected and moreover that there is a directed path from the root (the vertex corresponding to the unique continued fraction expansion) $\mathbf{b}_+$ to any other vertex in $G_{p/q}$.

\begin{figure}[h!]
    \centering
    \labellist \small
    \pinlabel 
    {{\color{my_green} 1+[-2,-3]}} at 90 45    \pinlabel{{\color{my_green} [2,-2]}} at 282 45
    \pinlabel{{\color{my_red} 2}} at 277 88
    \pinlabel{{\color{my_red} 2}} at 204 88
    \pinlabel {{\color{my_green} [3,2]}} at 351 45
    \pinlabel {complexity} at 292 230
    \pinlabel {slope} at 425 20
    \pinlabel {-4} at 157 -10
    \pinlabel {-2} at 199 -10
    \pinlabel {2} at 270 -10
    \pinlabel {4} at 308 -10
    \pinlabel {1} at 220 43
    \pinlabel {2} at 220 78
    \pinlabel {3} at 220 113
    \endlabellist
    \includegraphics[width=.47\textwidth]{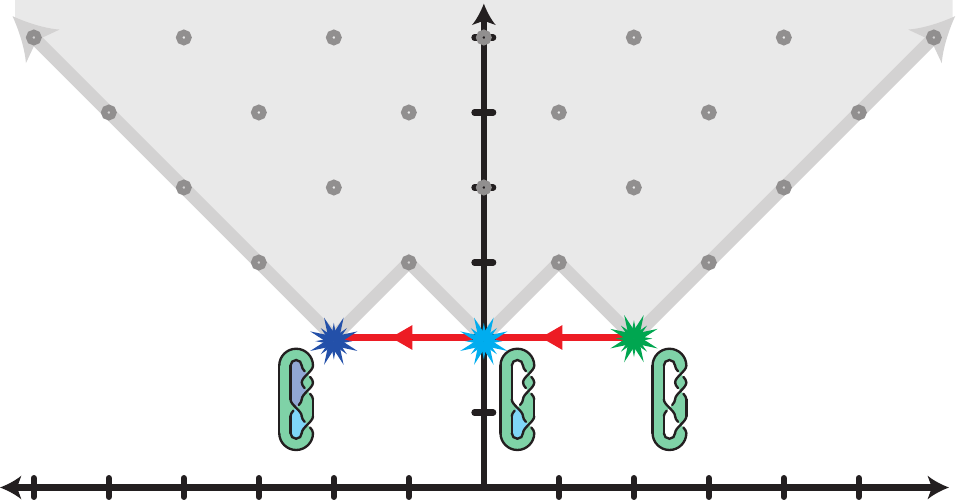}
\hspace{5mm}
    \labellist \small
    \pinlabel {complexity} at 208 230
    \pinlabel {slope} at 410 20
    \pinlabel {{\color{my_green} 1+[-2,-2]}} at 74 60
    \pinlabel {{\color{my_green} [3]}} at 301 30
    \pinlabel {{\color{my_orange} 3}} at 210 45
%    \pinlabel {-4} at 70 -10
%    \pinlabel {-2} at 113 -10
    \pinlabel {2} at 185 -10
    \pinlabel {4} at 221 -10
    \pinlabel {6} at 257 -10
    \pinlabel {1} at 164 43
    \pinlabel {2} at 164 78
    \pinlabel {3} at 164 113
    \endlabellist
    \includegraphics[width=.45\textwidth]{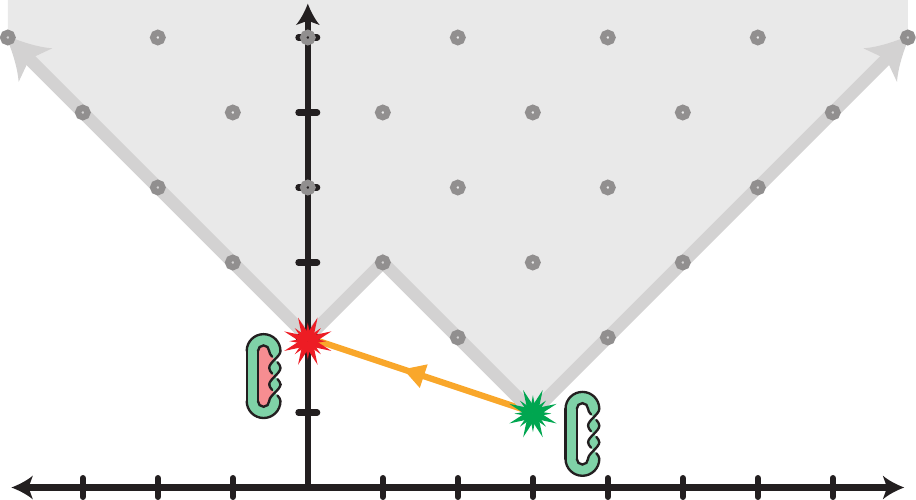}
    \caption{The geography  of the figure-eight knot $K_{2/5}$ on the left and that of the trefoil $K_{1/3}$ on the right}
    \label{Fi:Trefoil_Geography}
\end{figure}

\begin{figure}[h!]
    \centering
    \labellist\small
    \pinlabel{{\color{my_gray}{complexity}}} [l] at 152 395
    \pinlabel{{\color{my_gray} {$1$}}}  at 216 40
    \pinlabel{{\color{my_gray} {$2$}}}  at 216 76
    \pinlabel{{\color{my_gray} {$3$}}}  at 216 112
    \pinlabel{{\color{my_gray} {slope}}}  at 575 15
    \pinlabel{{\color{my_gray} {$-8$}}}  at 82 -8
    \pinlabel{{\color{my_gray} {$-6$}}}  at 118 -8
    \pinlabel{{\color{my_gray} {$-4$}}}  at 154 -8
    \pinlabel{{\color{my_gray} {$-2$}}}  at 190 -8
    \pinlabel{{\color{my_gray} {$4$}}}  at 298 -8
    \pinlabel{{\color{my_gray} {$8$}}}  at 370 -8  
    \pinlabel{{\color{my_gray} {$12$}}}  at 442 -8
    \pinlabel{{\color{my_gray} {$2$}}}  at 262 -8  
    \pinlabel{{\color{my_gray} {$6$}}}  at 334 -8
    \pinlabel{{\color{my_gray} {$10$}}}  at 406 -8  
    \pinlabel{{\color{my_gray} {$14$}}}  at 478 -8
    \tiny
    \pinlabel \rotatebox{90}{\scalebox{.8}{{\color{my_green} $\boldsymbol{1+[-2,-3,-4,-3,-2,-2]}$}}} at 83 150
    \pinlabel \rotatebox{90}{\scalebox{.8}{{\color{my_green} $\boldsymbol{[2,-2,-4,-3,-2,-2]}$}}} at 147 305
    \pinlabel \rotatebox{90}{\scalebox{.8}{{\color{my_green} $\boldsymbol{1+[-2,-3,-4,-2,3]}$}}} at 191 127
    \pinlabel \rotatebox{90}{\scalebox{.8}{{\color{my_green} $\boldsymbol{[3,2,-3,-3,-2,-2]}$}}} at 218 168
    \pinlabel \rotatebox{90}{\scalebox{.8}{{\color{my_green} $\boldsymbol{[2,-2,-4,-2,3],~1+[-2,-3,-3,2,4]}$}}} at 262 93
    \pinlabel \rotatebox{90}{\scalebox{.8}{{\color{my_purple} $\boldsymbol{1+[(-2)^{[2]},2^{[3]},(-2)^{[3]}%,-2,2,2,2,-2,-2,-2
    ]}$}}} at 238
    350
    \pinlabel \rotatebox{90}{\scalebox{.8}{{\color{my_green} $\boldsymbol{[3,2,-3,-2,3],~[2,-2,-3,2,4]}$}}} at 335 105
    \pinlabel \rotatebox{90}{\scalebox{.8}{{\color{my_green} $\boldsymbol{1+%[-2,-2,2,2,3,4]
    [(-2)^{[2]},2^{[2]},3,4]}$}}} at 378 283
    \pinlabel \rotatebox{90}{\scalebox{.8}{{\color{my_green} $\boldsymbol{%[3,3,2,2,-2,-2,-2]
    [3^{[2]},2^{[2]},(-2)^{[3]}]}$}}} at 327 310
    \pinlabel \rotatebox{90}{\scalebox{.8}{{\color{my_green} $\boldsymbol{[3,2,-2,2,4]}$}}} at 407 145
    \pinlabel \rotatebox{90}{\scalebox{.8}{{\color{my_green} $\boldsymbol{[3,3,2,3,4]}$}}} at 480 152
    \pinlabel {\scalebox{.8}{{\color{my_orange} $\boldsymbol{3}$}}} at 140 197
    \pinlabel {\scalebox{.8}{{\color{my_orange} $\boldsymbol{3}$}}} at 212 197
    \pinlabel {\scalebox{.8}{{\color{my_orange} $\boldsymbol{-3}$}}} at 271 243
    \pinlabel {\scalebox{.8}{{\color{my_orange} $\boldsymbol{3}$}}} at 277 283
    \pinlabel {\scalebox{.8}{{\color{my_orange} $\boldsymbol{3}$}}} at 277 210
    \pinlabel {\scalebox{.8}{{\color{my_orange} $\boldsymbol{-3}$}}} at 317 210
    \pinlabel {\scalebox{.8}{{\color{my_orange} $\boldsymbol{3}$}}} at 420 210
    \pinlabel {\scalebox{.8}{{\color{my_red} $\boldsymbol{2}$}}} at 232 180
    \pinlabel {\scalebox{.8}{{\color{my_red} $\boldsymbol{2}$}}} at 304 180
    \pinlabel {\scalebox{.8}{{\color{my_red} $\boldsymbol{2}$}}} at 376 180
    \pinlabel {\scalebox{.8}{{\color{my_red} $\boldsymbol{2}$}}} at 448 180
    \pinlabel {\scalebox{.8}{{\color{my_red} $\boldsymbol{2}$}}} at 114 216
    \pinlabel {\scalebox{.8}{{\color{my_red} $\boldsymbol{2}$}}} at 186 216
    \pinlabel {\scalebox{.8}{{\color{my_cyan} $\boldsymbol{-4}$}}} at 164 257
    \pinlabel {\scalebox{.8}{{\color{my_cyan} $\boldsymbol{4}$}}} at 290 255
    \pinlabel {\scalebox{.8}{{\color{my_cyan} $\boldsymbol{4}$}}} at 415 225

    \endlabellist
\includegraphics[width=\textwidth]{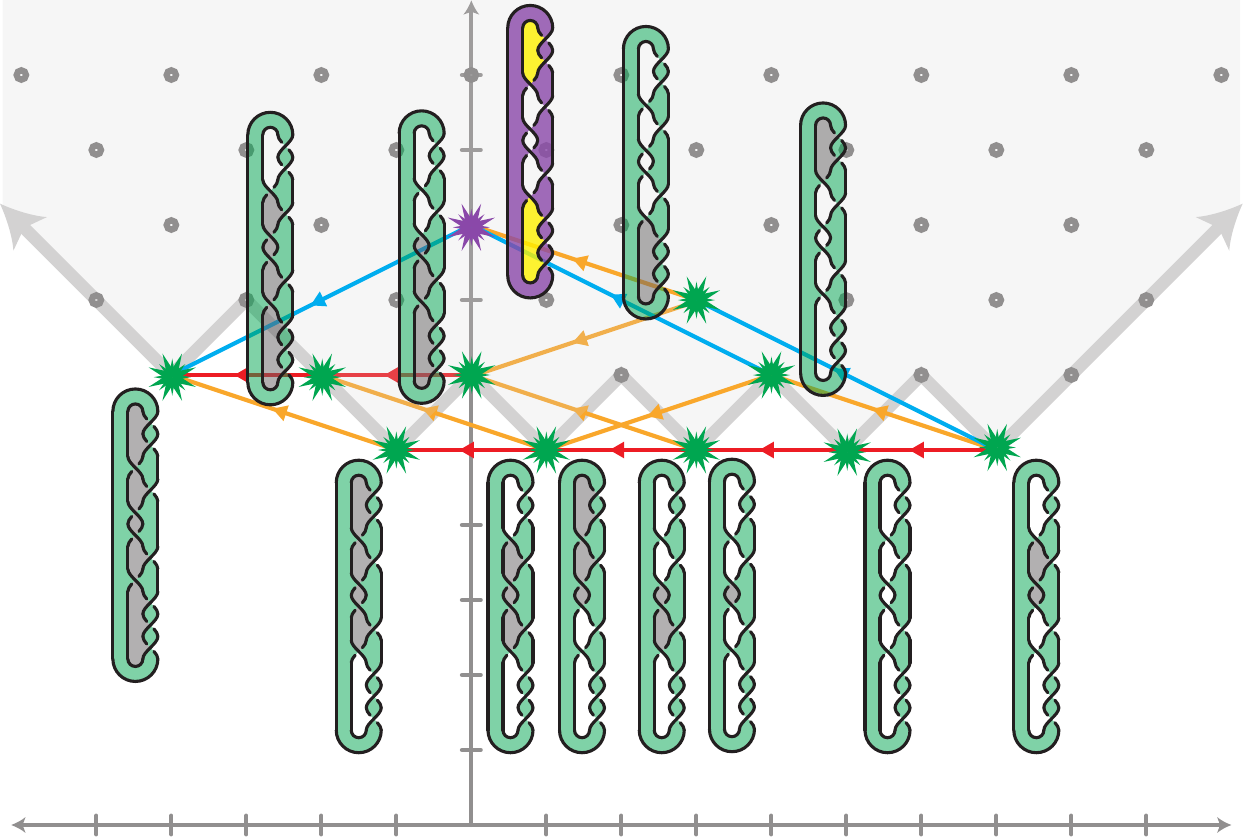}
    \caption{The image of the graph $G_{43/111}$ projected onto the geography of the knot $K_{43/111}$ }
    \label{Fi:43/111_Geography}
\end{figure}

Before closing the introduction with a brief outline, we provide three motivating examples, the last of which will continue as a running example throughout the paper. In Figures~\ref{Fi:Trefoil_Geography} and \ref{Fi:43/111_Geography}, we depict the set of all pairs $(\beta_1,s)$ where $\beta_1$ is the first Betti number and $s$ is the boundary slope of a spanning surface $F$ for the given knot. This set of pairs is called the {\it geography} of the knot (or link). For example, see \cite{Allen, kink,Howie, Sabloff}.   
Each starred vertex in these figures represents a complexity-slope pair corresponding to a $\pi_1$-essential spanning surface.

There is a natural map from the graph $G_{p/q}$ to the geography for $K_{p/q}$ that sends each vertex of $G_{p/q}$ corresponding to an essential continued fraction expansion $\mathbf{b}$ of $\frac{p}{q}$ to the starred vertex $(\beta_1(F_\mathbf{b}),\text{slope}(F_{\mathbf{b}}))$, where $F_{\mathbf{b}}$ is an essential spanning surface for $K_{p/q}$ associated to $\mathbf{b}$.  This map allows one to see $G_{p/q}$ (or a quotient of it) within the geography picture for $K_{p/q}$. Again, see Figures~\ref{Fi:Trefoil_Geography} and \ref{Fi:43/111_Geography}.

\begin{example}\label{Ex:2/5}
As a warm-up, consider the figure-eight knot $K_{2/5}$. 
There are three continued fraction expansions of $\frac{2}{5} = r + [b_1, \ldots, b_k]$ that are essential (again, meaning that $|b_i|\geq 2$ for all $i$), namely 
\[1+[-2,-3]:=1+\cfrac{1}{-2-\cfrac{1}{-3}},\qquad [2,-2]:=\cfrac{1}{2-\cfrac{1}{-2}},\qquad\text{ and }\qquad[3,2]:=\cfrac{1}{3-\cfrac{1}{2}}.\]
In this specific instance, all three essential continued fractions are related by Hirasawa-Teragaito's 2-moves, which we define in \textsection\ref{subsec:HT}. The two red edges labeled by 2 shown on the left in Figure \ref{Fi:Trefoil_Geography} depict the corresponding relationships among the three essential spanning surfaces for $K_{2/5}$.  These two edges and the three vertices are the graph $G_{\frac{2}{5}}$.
\end{example}

\begin{example}\label{Ex:1/3}
In the even simpler case of the  trefoil $K_{1/3}$, 
there are just two essential expansions of $\frac{1}{3}$, namely 
\[[3]:= \frac{1}{3}\qquad\text{ and }\qquad1+[-2,-2]:= 1+\cfrac{1}{-2-\cfrac{1}{-2}},\] 
and the only ways to relate these using just Hirasawa-Teragaito's moves involve inessential continued fractions. The expansions $[3]$ and $1+[-2,-2]$ are directly related to each other, however, by a $3$-move, which we will define in \textsection \ref{S:n-moves}. The orange edge labeled by 3 shown on the right in Figure \ref{Fi:Trefoil_Geography} depicts the corresponding relationship between the two essential spanning surfaces for $K_{1/3}$.  This edge and the two vertices are the graph $G_{\frac{1}{3}}$.
\end{example}

We now introduce the paper's main running example.

\begin{example}\label{Ex:43/111_intro}
Figure \ref{Fi:43/111_Geography} shows the more complicated example of $K_{43/111}$ which appears in the knot tables as $11a176$. This knot has 13 essential continued fraction expansions which are mutually related by $\pm 2$-moves (red), $\pm3$-moves (orange), and $\pm4$-moves (blue). We will expand upon this in Examples \ref{Ex:Euclidean},  
\ref{Ex:43/111_Hopf}, \ref{Ex:43/111_moves}, and  \ref{Ex:43/111_G_p/q}.
\end{example}

Here is a brief outline of the paper. 
In \textsection\ref{S:background1}, we review continued fraction expansions of rational numbers and Hirasawa-Teragaito's moves on them.
In \textsection\ref{S:moves}, we define $\pm n$-moves on continued fraction expansions and prove Theorem \ref{T:essentialconnected}.
In \textsection\ref{S:Examples}, we collect the graphs $G_{p/q}$ for all knots and links $K_{p/q}$ with seven or fewer crossings, characterize the graphs $G_{p/q}$ precisely for several infinite families, and examine a few properties of the graph $G_{p/q}$.

In Appendix~\ref{allmoves}, we list explicit formulas for Hirasawa-Teragaito's $0$-, $\pm1$-, $232$-, and $\pm2$-moves and our $\pm n$-moves, including all edge cases (in the main text, we provide only one formula for each move, since the edge cases are similar, although not identical).

\subsection*{AI Declaration}  
No AI tools were used in this project.

\section{Background on continued fractions}\label{S:background1}

\subsection{Subtractive continued fraction expansions of rational numbers}

Every sequence of integers $r,b_1,b_2,\dots, b_k$ gives rise to a reduced fraction $\frac{p}{q}$ via a {\bf subtractive continued fraction}
\begin{equation}\label{E:CF1}
    \frac{p}{q}= r+[b_1, b_2,\dots, b_k] := r + \cfrac{1}{b_1-\cfrac{1}{b_2-\cfrac{1}{\ddots -\cfrac{1}{b_k}}}}.
\end{equation}
Sometimes in this computation (e.g. if $b_k=0$ or $b_k=1=b_{k-1}$), one encounters an expression of the form 
$\cfrac{1}{b-\frac{1}{0}}$. This should be treated as equaling 0.  We also declare that $r+[0]=\frac{1}{0}$ for any integer $r$. Rational numbers $\frac{p}{q}$ will always be written in lowest terms with $q>0$, and continued fractions will always be subtractive, denoted as in Equation \eqref{E:CF1}. 

Conversely, every rational $\frac{p}{q}\in\mathbb{Q}\setminus \mathbb{Z}$ has a continued fraction expansion as in Equation~\eqref{E:CF1}. The continued fraction expansion for $\frac{p}{q}$ is not unique. There are two particular continued fraction expansions of note associated to each $\frac{p}{q}\in\mathbb{Q}\setminus \mathbb{Z}$, namely the positive continued fraction expansion and the negative continued fraction expansion.

\begin{definition}
\label{def:posneg}
A continued fraction expansion $\frac{p}{q} =r+[b_1,\dots,b_k]$ is called \textbf{positive}, (respectively,  \textbf{negative}) if $b_i\geq 2$ (respectively, $b_i \le -2$) for all $i$.
\end{definition}

Every rational $\frac{p}{q}\in\mathbb{Q}\setminus\mathbb{Z}$ has both a positive and a negative continued fraction expansion, and those expansions are unique. We denote them $\mathbf{b}_+$ and $\mathbf{b}_-$, respectively. We do not include the proof of uniqueness here, but the proof runs parallel to the classical proof of the same result for additive continued fractions (see, for example, Theorem 160 in~\cite{hardy}).

To find the negative continued fraction expansion $\mathbf{b}_-$ for $\frac{p}{q}$, apply a modified version of the Euclidean algorithm, which uses the ceiling function instead of the floor function. 
To find a positive continued fraction expansion $\mathbf{b}_+$ for $\frac{p}{q}$, use the modified Euclidean algorithm starting with the second step. These processes are illustrated in the following example. 

\begin{example}
    \label{Ex:Euclidean}
Consider $\frac{43}{111}$ from Example \ref{Ex:43/111_intro}.
    In the left (respectively, right) column, we use a modified Euclidean algorithm with the ceiling function to find the positive (respectively, negative) continued fraction expansion. 
{\small
\begin{equation*}
\begin{split}
43 & = (0)111 +43\\
111 & =(3)43 -18\\
43 & = (3)18 - 11 \\
18 & = (2)11 -4 \\
11 &= (3)4 -1 \\
4 &= (4)1 -0 \\
&~~
\end{split}
\quad
\quad
\quad
\begin{split}
43 & = (1)111 -68 \\
111 & =(2)68 -25\\
68 & = (3)25 - 7 \\
25 & = (4)7 -3 \\
7 &= (3)3 -2 \\
3 &= (2)2 -1\\
2 &= (2)1 -0
\end{split}
\end{equation*}
}
Therefore,  $\mathbf{b}_+ = 0+[3,3,2,3,4]$ and $\mathbf{b}_- = 1+[-2,-3,-4,-3,-2,-2]$. 
\end{example}

It will be convenient to introduce the notation $b^{[m]}$ to represent a string of $m\geq 0$ consecutive $b$'s and where $b^{[0]}$ represents an empty string. Thus, in Example \ref{Ex:Euclidean}, we could also write $\frac{43}{111} = 0+[3^{[2]},2,3,4]$ or even $\frac{43}{111} = 0+[3^{[2]},2,(-2)^{[0]},3,4]$. 

\subsection{Hirasawa and Teragaito's moves on continued fractions}
\label{subsec:HT}

Motivated by a result of Bessho \cite{Bessho_1994} about nonorientable spanning surfaces for 2-bridge knots and links, Hirasawa and Teragaito describe three moves on a subtractive continued fraction $r+[b_1,\dots,b_k]$ that shorten its length without changing the value of the associated rational number \cite{Hirasawa_Teragaito:2006}. We name these moves for ease of reference.

    \begin{itemize}[leftmargin=*, widest=\textbf{232-move: }]
        \item[\textbf{0-move: }] Removal of term 0   
        \begin{equation*}\label{E:0}
            r+[\ldots,a,0,b,\ldots] \rightarrow r+[\ldots,a+b,\ldots]
        \end{equation*}
        \item[\textbf{$\boldsymbol{\pm 1}$-move: }] Removal of term  $\upvarepsilon=\pm1$
        \begin{equation*}\label{E:1}
            r+[\ldots,a,\upvarepsilon,b,\ldots] \rightarrow r+[\ldots,a-\upvarepsilon,b-\upvarepsilon,\ldots]
        \end{equation*}
        \item[\textbf{232-move: }] Replacement of a subsequence of the form $2\upvarepsilon, (3\upvarepsilon)^{[m-2]}, 2\upvarepsilon$ for $m\geq2$, $\upvarepsilon=\pm1$
        \begin{equation*}\label{E:232}
            r+[\ldots,a,2\upvarepsilon, (3\upvarepsilon)^{[m-2]}, 2\upvarepsilon,b,\ldots] \rightarrow r+[\ldots,a-\upvarepsilon,(-3\upvarepsilon)^{[m-1]}, b-\upvarepsilon,\ldots]
        \end{equation*}
    \end{itemize}

Each formula above governs the cases where the term or subsequence being removed or replaced is neither at the very start nor the very end of the continued fraction expansion.  Similar formulas hold in the edge cases. We provide explicit formulas for all of these cases in Appendix~\ref{allmoves}.

Hirasawa and Teragaito use the three moves defined above to obtain the following result.

\begin{theorem}[Hirasawa-Teragaito, Theorem 2] \label{Alg:HT}
    Let $\frac{p}{q}=r+[b_1,\ldots,b_k]$ be a continued fraction expansion obtained from an arbitrary expansion of $\frac{p}{q}$ by fully reducing the length by a (finite) repetition of the 0-, $\pm$1-, and 232-moves above.  
    Then $k$ is the minimal length among all continued fraction expansions of $\frac{p}{q}$.
\end{theorem}

Theorem \ref{Alg:HT} provides an avenue for finding {\it some} continued fraction expansion of minimal length for a given rational number $\frac{p}{q}$. Hirasawa-Teragaito showed further that {\it all} minimal length continued fraction expansions of $\frac{p}{q}$ can be obtained from this particular one by repeatedly performing an additional pair of moves, which we call {\bf $\boldsymbol{\pm2}$-moves}.
We again write these moves just in the cases where the $\pm2$ occurs somewhere other than the first or last term. Explicit formulas for boundary cases are provided in Appendix~\ref{allmoves}.

\begin{itemize}[leftmargin=*, widest=\textbf{$\boldsymbol{\pm2}$-move: }]        
\item[\textbf{$\boldsymbol{\pm2}$-move: }] Replacing a term $2 \varepsilon $ with $-2 \varepsilon $ where $\varepsilon = \pm 1$
        \begin{equation*}
            r+[\ldots,a,2 \varepsilon ,b,\ldots] \to r+[\ldots,a-\varepsilon,-2\varepsilon, b-\varepsilon,\ldots]
        \end{equation*}
\end{itemize}

\begin{theorem}[Hirasawa-Teragaito, Theorem 4] \label{T:HT2}
    Any two shortest continued fraction expansions for any $\frac{p}{q}\in\mathbb{Q}\setminus\mathbb{Z}$ are related to each other by a finite sequence of $\pm2$-moves.
\end{theorem}

\begin{example}
\label{Ex:43/111_Hopf}
In our running Example \ref{Ex:43/111_intro}, starting from $\mathbf{b_+}=[3,3,2,3,4]$, one can iteratively perform six $\pm2$-moves:
\begin{align*}    [3,3,\boxed{2},3,4]&\underset{2\text{-move}}{\longrightarrow}[3,\boxed{2},-2,2,4]\underset{2\text{-move}}{\longrightarrow}[\boxed{2},-2,-3,2,4]
    \underset{2\text{-move}}{\longrightarrow} 1+[-2,-3,-3,\boxed{2},4]\\
    &\underset{2\text{-move}}{\longrightarrow}1+[\boxed{-2},-3,-4,-2,3]
    \underset{-2\text{-move}}{\longrightarrow}[2,\boxed{-2},-4,{-2},3]
    \underset{-2\text{-move}}{\longrightarrow}[3,2,-3,-2,3].
\end{align*}
These are, in fact, all minimal-length continued fraction expansions of $\frac{43}{111}$, so we have verified Theorem \ref{T:HT2} in this example.  See the red path along the bottom of Figure \ref{Fi:43/111_Geography}. 
\end{example}

\section{Moves connecting all essential continued fraction expansions of a rational number}
\label{S:moves}

\subsection{Essential continued fraction expansions}
 
The following definition is, as noted in the introduction, motivated by part (c) of Hatcher-Thurston's Theorem 1 in \cite{Hatcher-Thurston:1985}, which states which continued fraction expansions correspond to $\pi_1$-essential surfaces.

\begin{definition}
\label{def:essential}
We call a subtractive continued fraction $r+[b_1,\dots,b_k]$ {\bf essential} if $|b_i|\geq 2$ for each $i=1,\dots,k$.    Otherwise such a continued fraction is \textbf{inessential}.
\end{definition}

\begin{remark}\label{R:integer_part}
One convenient property of an essential continued fraction expansion $\frac{p}{q}=r+[b_1,\dots,b_k]$ is that the sign of the first term $b_1$ determines the integer $r$. Namely, $r=\lfloor\frac{p}{q}\rfloor$ if $b_1>0$, and $r=\lceil \frac{p}{q}\rceil$ if $b_1<0$. In particular, when $0<\frac{p}{q}<1$, we have $r=0$ when $b_1>0$ and $r=1$ when $b_1<0$. 
\end{remark}

\subsection{Moves on essential continued fraction expansions}\label{S:n-moves}

We now generalize Hirasawa and Teragaito's $\pm2$-moves.

\begin{definition}\label{def:n_move} 
Consider an essential continued fraction expansion 
$\mathbf{b}=r+[b_1,\dots,b_k]$, and let $n$ be an integer greater than or equal to 2. An {\bf $\boldsymbol{n}$-move} (on the $i^\text{th}$ term of $\mathbf{b}$) replaces $b_i=n$ with a string of the form $(-2)^{[n-1]}$, and a {\bf $\boldsymbol{-n}$-move} replaces $b_i=-n$ with a string of the form $2^{[n-1]}$; in each case, one must also adjust the $(i-1)^\text{st}$ and $(i+1)^\text{st}$ terms.  Here are the formulas when $i=2,\dots,k-1$.

\begin{itemize}[leftmargin=*, widest=\textbf{$-n$-move: }]
        \item[\textbf{$\boldsymbol{\pm n}$-move: }] Replacement of term $b_i=n \varepsilon$ 
        with subsequence $(-2 \varepsilon)^{[n-1]}$ for $n\geq 2$, $\varepsilon = \pm 1$  
        \begin{equation*}
            r+[b_1,\dots,b_{i-1},n\varepsilon,b_{i+1},\dots,b_k]\to r+[b_1,\dots,b_{i-1}-\varepsilon,(-2 \varepsilon)^{[n-1]},b_{i+1}-\varepsilon,\dots,b_k]
        \end{equation*}
\end{itemize}

When $i=1,k$, the same formulas hold, but with the first or last term truncated, respectively; if the first term is truncated, then $r$ also needs to be adjusted as specified by Remark \ref{R:integer_part}. We provide explicit formulas for all cases in Appendix~\ref{allmoves}.

Inverse $\pm n$-moves thus have the following form:
\begin{equation}\label{E:inverse_move}
                r+[b_1,\dots,b_{i},(-2\varepsilon)^{[n-1]},b_{j},\dots,b_k] \overset{\text{inverse}}{\underset{n \varepsilon \text{-move}}{\longrightarrow}}~ r+[b_1,\dots,b_i+\varepsilon,n\varepsilon,b_j+\varepsilon,\dots,b_k].
\end{equation}

Call an inverse $n\varepsilon$-move as in Equation \eqref{E:inverse_move} {\bf maximal} if $b_i,b_j\neq -2\varepsilon$. Call a $\pm n$-move or inverse $\pm n$-move $\mathbf{b}\to\mathbf{b}'$ {\bf essential} if $\mathbf{b}$ and $\mathbf{b}'$ are both essential.
\end{definition}

\begin{example}\label{Ex:43/111_moves}  We continue our running example that began in Example \ref{Ex:43/111_intro} 
 by considering the five different $n$-moves that one could perform on the expansion $\mathbf{b_+}=[3,3,2,3,4]$ of $\frac{p}{q}=\frac{43}{111}$. Three of these moves are essential, 
while each of the other two is inessential, 
as indicated below by a \textcircled{1} on the right.
\begin{align*}
[\text{\boxed{3}},3,2,3,4]\underset{3\text{-move}}{\longrightarrow}&1+[-2,-2,2,2,3,4],\\
    [3,\text{\boxed{3}},2,3,4]\underset{3\text{-move}}{\longrightarrow}&[2,-2,-2,\text{\textcircled{1}},3,4],\\
    [3,3,\text{\boxed{2}},3,4]\underset{2\text{-move}}{\longrightarrow}&[3,2,-2,2,4],\\
    [3,3,2,\text{\boxed{3}},4]\underset{3\text{-move}}{\longrightarrow} &[3,3,\text{\textcircled{1}},-2,-2,3],\\
    [3,3,2,3,\text{\boxed{4}}]\underset{4\text{-move}}{\longrightarrow}&[3,3,2,2,-2,-2,-2].
\end{align*}
\end{example}

Given an essential continued fraction expansion $\mathbf{b}$ of a rational number $\frac{p}{q}$, it is always possible to find some essential $n$-move or inverse $n$-move $\mathbf{b}\to\mathbf{b}'$.

\begin{proposition}
\label{P:essential}
    Suppose $\mathbf{b}=r+[b_1,\dots,b_k]$ is an essential continued fraction expansion of a non-integer rational number. 
    \begin{enumerate}[label=(\arabic*)]
        \item\label{part:essential_1} Given an inverse $\pm n$-move $\mathbf{b}\to\mathbf{b}'$ for any $n\geq 2$, $\mathbf{b}'$  is essential if and only if this move is maximal.
        \item\label{part:essential_2_and_3} Given an $\varepsilon n$-move $\mathbf{b}\to\mathbf{b}'$ on the $i^{\text{th}}$ term of $\mathbf{b}$ for any $n\geq 2$ and $\varepsilon=\pm1$, the continued fraction $\mathbf{b}'$ is essential if and only if $b_{i-1}\neq2\varepsilon$ (or $i=1$) and $b_{i+1}\neq 2\varepsilon$ (or $i=k$).
        \item\label{part:essential_4-} 
        If $\mathbf{b}\neq\mathbf{b}_+$, then $\mathbf{b}$ admits an essential $-n$-move or essential inverse $n$-move for some $n$.
        \item\label{part:essential_4+} 
        If $\mathbf{b}\neq\mathbf{b}_-$, then $\mathbf{b}$ admits an essential $n$-move or essential inverse $-n$-move for some $n$.
    \end{enumerate}
\end{proposition}

\begin{proof}
    For  part \ref{part:essential_1}, the result of Equation \eqref{E:inverse_move} is inessential if and only if $b_i=-2\varepsilon$ or $b_j=-2\varepsilon$, which is precisely when the move is non-maximal. Part \ref{part:essential_2_and_3} follows immediately from Definitions \ref{def:essential} and \ref{def:n_move}. 

    For part \ref{part:essential_4-}, if some term of $\mathbf{b}$ equals $-2$, then choose the largest $n$ for which $\mathbf{b}$ admits an inverse $n$-move; part \ref{part:essential_1} ensures that the resulting continued fraction is essential.  
    Otherwise, choose a negative term $b_i=-n$ and perform an $-n$-move on this $i^{\text{th}}$ term of $\mathbf{b}$; since every negative term of $\mathbf{b}$ is at most $-3$, part \ref{part:essential_2_and_3} ensures that the resulting continued fraction is essential. The proof of part \ref{part:essential_4+} is almost identical.
\end{proof}

\begin{definition}\label{D:G_p/q}
Given $\frac{p}{q}\in\mathbb{Q}\setminus\mathbb{Z}$, let $G_{p/q}$ denote the directed graph whose vertices correspond to essential continued fraction expansions of $\frac{p}{q}$ and whose edges correspond to essential $\pm n$-moves between these expansions. Orient each edge forward for an $n$-move and backward for a $-n$-move. 
\end{definition}

\begin{figure}[h!]
\centering
\labellist\tiny
\pinlabel \rotatebox{26.5}{ {\color{my_green} $\mathbf{b}_-=1+[-2,-3,-4,-3,-2,-2]$} } at 30 64
\pinlabel{ {\color{my_green} $[2,-2,-4,-3,-2,-2]$} } at 85 50
\pinlabel{ {\color{my_green} $[3,2,-3,-3,-2,-2]$} } at 195 45
\pinlabel{ {\color{my_purple} $1+[-2,-2,2,2,2,-2,-2,-2]$} } at 150 122
\pinlabel %\rotatebox{-26.5}
{ {\color{my_green} $[3,3,2,2,-2,-2,-2]$} } at 285 87
\pinlabel \rotatebox{-20}{ {\color{my_green} $1+[-2,-2,2,2,3,4]$} } at 300 47
\pinlabel{ {\color{my_green} $1+[-2,-3,-4,-2,3]$} } at 75 5
\pinlabel{ {\color{my_green} $1+[-2,-3,-3,2,4]$ }} at 150 -5
\pinlabel{ {\color{my_green} $[2,-2,-3,2,4]$ }} at 250 -5
\pinlabel\rotatebox{-17}{{ \color{my_green} $[2,-2,-4,-2,3]$} } at 165 20
\pinlabel{ {\color{my_green} $[3,2,-3,-2,3]$} } at 280 13
\pinlabel{ {\color{my_green} $[3,2,-2,2,4]$} } at 330 0
\pinlabel{ {\color{my_green} $\mathbf{b}_+=[3,3,2,3,4]$} \color{black}} at 400 0
\pinlabel{{\color{my_cyan} $-4$}} at 78 77
\pinlabel{{\color{my_cyan} $4$}} at 221 75
\pinlabel{{\color{my_cyan} $4$}} at 332 45
\pinlabel{{\color{my_orange} $3$}} at 58 21
\pinlabel{{\color{my_orange} $3$}} at 132 18
\pinlabel{{\color{my_orange} $3$}} at 203 18
\pinlabel{{\color{my_orange} $-3$}} at 205 59
\pinlabel{{\color{my_orange} $3$}} at 208 100
\pinlabel{{\color{my_orange} $-3$}} at 239 20
\pinlabel{{\color{my_orange} $3$}} at 346 21
\pinlabel{{\color{my_red} $2$}} at 43 40
\pinlabel{{\color{my_red} $2$}} at 115 40
\pinlabel{{\color{my_red} $2$}} at 150 1
\pinlabel{{\color{my_red} $2$}} at 222 -2
\pinlabel{{\color{my_red} $2$}} at 296 1
\pinlabel{{\color{my_red} $2$}} at 366 12
\endlabellist
\includegraphics[width=\textwidth]{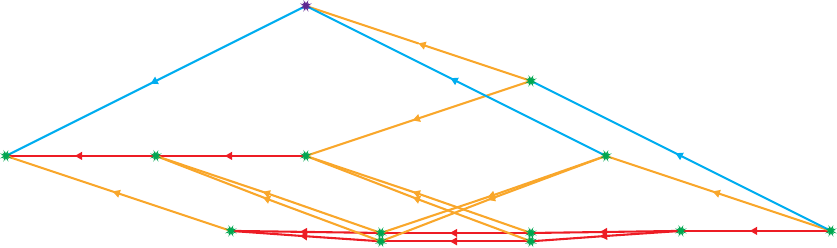}
\caption{The graph $G_{43/111}$}
\label{Fi:G_43/111}
\end{figure}

Parts \ref{part:essential_4-} and \ref{part:essential_4+} of Proposition \ref{P:essential} imply that every vertex in $G_{p/q}$ has positive degree. This is because the degree of any vertex in $G_{p/q}$ will equal the number of essential $\pm n$-moves and inverse $\pm n$-moves on the corresponding continued fraction expansion. In Example \ref{Ex:43/111_moves}, e.g., we saw that $\mathbf{b}_+=[3,3,2,3,4]$ admits three essential $n$-moves (and no $-n$-moves), which is reflected in $G_{43/111}$ (see Figure \ref{Fi:G_43/111}), where the vertex  corresponding to $\mathbf{b}_+$ has degree three.

Our main result, Theorem \ref{T:essentialconnected}, will imply that $G_{p/q}$ is always connected, and moreover that, given any vertex $v$ in $G_{p/q}$, there is a directed path in $G_{p/q}$ that begins at the vertex corresponding to $\mathbf{b}_+$ and ends at $v$. First, we build further upon our running example. 

\begin{example}\label{Ex:43/111_G_p/q}
Figure \ref{Fi:G_43/111} shows the graph $G_{43/111}$ (c.f. Figure \ref{Fi:43/111_Geography}). There are many paths in $G_{43/111}$ between  $\mathbf{b_+}=[3,3,2,3,4]$ and $\mathbf{b_-}=1+[-2,-3,-4,-3,-2,-2]$, including
\begin{align*}
    [3,3,\boxed{2},3,4]&\underset{2\text{-move}}{\longrightarrow}[3,\boxed{2},-2,2,4]\underset{2\text{-move}}{\longrightarrow}[\boxed{2},-2,-3,2,4]
    \underset{2\text{-move}}{\longrightarrow} 1+[-2,-3,-3,\boxed{2},4]\\
    &\underset{2\text{-move}}{\longrightarrow}1+[-2,-3,-4,-2,\boxed{3}]
    \underset{3\text{-move}}{\longrightarrow}1+[-2,-3,-4,-3,-2,-2]
\end{align*}
 along the very bottom of the graph and 
\begin{align*}
    [3,3,2,3,\boxed{4}]&\underset{4\text{-move}}{\longrightarrow}[\boxed{3},3,2,2,-2,-2,-2]\underset{3\text{-move}}{\longrightarrow}1+[-2,-2,\boxed{2,2,2},-2,-2,-2]\\
    &\underset{\text{inverse }(-4)\text{-move}}{\longrightarrow}1+[-2,-3,-4,-3,-2,-2]
\end{align*}
along the very top.  Note that both of these paths are directed, proceeding from right to left.
\end{example}

 \subsection{Main result}
 
Our main theorem, as stated in the introduction, assures that the types of paths described at the end of Example \ref{Ex:43/111_moves} always exist---that is, $G_{p/q}$ is connected as a directed graph with root vertex $\mathbf{b}_+$.

We are ready to prove our main result.

\begin{T:essentialconnected}
Let $\frac{p}{q}\in\mathbb{Q}\setminus\mathbb{Z}$. 
    If $\mathbf{b}_+=r+[b_1,\dots,b_k]$ is the unique continued fraction expansion of $\frac{p}{q}$ in which each term $b_i\geq 2$, then any essential continued fraction expansion of $\frac{p}{q}$ can be obtained from $\mathbf{b}_+$ via a sequence of essential $n$-moves and essential inverse $-n$-moves for various $n\geq 2$.  
\end{T:essentialconnected}

\begin{proof}
Part \ref{part:essential_4-} of Proposition \ref{P:essential} states that any continued fraction $\mathbf{b}$ other than $\mathbf{b}_+$ of $p/q$ admits an essential $(-n)$-move or an essential inverse $n$-move. An essential $(-n)$-move reduces the number of negative entries in $\mathbf{b}$ by one, and an essential inverse $n$-move reduces the number of negative entries in $\mathbf{b}$ by $n-1$. Thus there is a sequence of such moves transforming $\mathbf{b}$ into $\mathbf{b}_+$. Reversing this sequence of moves results in the desired sequence from $\mathbf{b}_+$ to $\mathbf{b}$.
\end{proof}

As noted in the introduction, Theorem \ref{T:essentialconnected} has the following immediate corollary:

\begin{C:essentialconnected}
        Any two essential continued fraction expansions of $\frac{p}{q}$ are related by a sequence of essential $\pm n$-moves and their inverses for various $n\geq 2$. 
\end{C:essentialconnected}

If $\frac{p}{q}=r+[b_1,b_2,\dots,b_k]$, then $[-b_1,-b_2,\dots,-b_k]\equiv-\frac{p}{q}\pmod{1}$.  Moreover, any $n$-move on $[b_1,\dots,b_k]$ corresponds to a $(-n)$-move on $[-b_1,\dots,-b_k]$, and likewise for inverse moves.  Theorem \ref{T:essentialconnected} thus implies the following.

\begin{corollary}\label{C:Essential}
    Let $\frac{p}{q}\in\mathbb{Q}\setminus\mathbb{Z}$. For any essential continued fraction expansion $\mathbf{b}'$ of $\frac{p}{q}$, there is a sequence of essential $n$-moves and essential inverse $-n$-moves taking $\mathbf{b}'$ to $\mathbf{b}_-$.
    Hence, there is a sequence of these moves taking $\mathbf{b}_+$ to $\mathbf{b}_-$ by way of $\mathbf{b}'$.
\end{corollary}

\section{Examples of graphs for 2-bridge knots and links}\label{S:Examples}

Next, we provide tables with many examples of the graphs $G_{p/q}$.  Those tables, together with the subsequent Remarks \ref{R:Negate} and \ref{R:Invert}, describe the graphs $G_{p/q}$ for all knots and links $K_{p/q}$ with seven or fewer crossings. After that, we highlight several infinite families for which the graph $G_{p/q}$ is simple to characterize, and we highlight an interesting behavior exhibited by the graphs $G_{q/(q^2+1)}$.

\subsection{Graphs for 2-bridge knots and links with small crossing number}
\label{subsec:crossing}
In Tables \ref{Table:0}  
and \ref{Table:2} in  this subsection, we collect the graphs $G_{p/q}$ for all equivalence classes of rational numbers $\frac{p}{q}$ corresponding to 2-bridge knots and links with seven crossings or fewer. The equivalence here is generated by the rule that $\frac{p}{q}\sim \frac{p'}{q}$ whenever $p$ and $q$ are in lowest terms, $q$ is positive, and either $p\equiv \pm p'\pmod{q}$ or $p p'\equiv \pm 1\pmod{q}$. That is, two rational numbers are identified if the corresponding 2-bridge knots or links are either equivalent or mirror images. If $r+[b_1,\dots,b_k]$ is the positive continued fraction expansion for $\frac{p}{q}$ then the crossing number of the associated knot or link $K_{p/q}$ is $b_1+\cdots+b_k+1-k$.  
See Cromwell \cite[p.251]{Cromwell2004}, Koseleff and Pecker \cite[Proposition 2.5]{KosPec}, Suzuki \cite[Theorem 3.3]{Suzuki}, or the preceding paper by the authors \cite[Section 2.3]{CohKinLowShaVan}.  
Horizontal dividing lines in Table \ref{Table:0} separate the examples by crossing number:  2, 3, 4, 5, and 6.
\begin{table}[h!]
\centering
\begin{tabular}{|| c | c  c  c ||}
 \hline
 \hline
2 & \labellist\tiny
\pinlabel{{\color{my_red} $2$}} at 59 11
\pinlabel{{\color{my_purple} $[-2]$}} at 19 26
\pinlabel{$\frac{1}{2}$} at 5 10
\pinlabel{{\color{my_purple} $[2]$}} at 89 26
\endlabellist
\includegraphics[scale=.8]{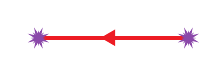}&&\\ 
\hline
3 & \labellist\tiny
\pinlabel \rotatebox{-45}{{\color{my_purple} $[(-2)^{[2]}]$}} at 28 22
\pinlabel{{\color{my_green} $[3]$}} at 118 24
\pinlabel{{\color{my_orange} $3$}} at 70 34
\pinlabel{$\frac{1}{3}$} at 5 10
\endlabellist
\includegraphics[scale=.7]{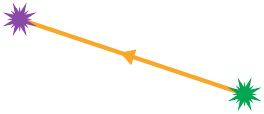}&&\\ \hline
4 & \labellist\tiny
\pinlabel \rotatebox{-45}{{\color{my_purple} $[(-2)^{[3]}]$}} at 27 52
\pinlabel{{\color{my_purple} $[4]$}} at 154 25
\pinlabel{$\frac{1}{4}$} at 7 12
\pinlabel{{\color{my_green} $4$}} at 90 49
\endlabellist
\includegraphics[scale=.6]{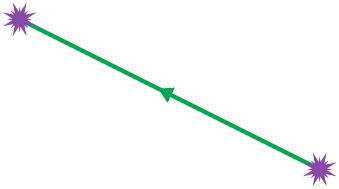}&
\labellist\tiny
\pinlabel{{\color{my_green} $[-2,-3]$}} at 30 37
\pinlabel{{\color{my_purple} $[2,-2]$}} at 85 37
\pinlabel{{\color{my_green} $[3,2]$}} at 145 37
\pinlabel{$\frac{2}{5}$} at 7 12
\pinlabel{{\color{my_red} $2$}} at 50 62
\pinlabel{{\color{my_red} $2$}} at 125 62
\endlabellist
\includegraphics[scale=.6]{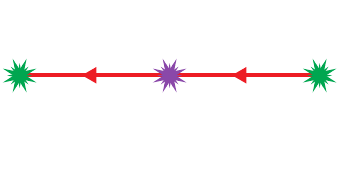}&\\ \hline
5 & \labellist\tiny
\pinlabel \rotatebox{-45}{{\color{my_purple} $[(-2)^{[4]}]$}} at 33 92
\pinlabel{{\color{my_green} $[5]$}} at 188 37
\pinlabel{$\frac{1}{5}$} at 7 14
\pinlabel{{\color{my_cyan} $5$}} at 110 75
\endlabellist
\includegraphics[scale=.5]{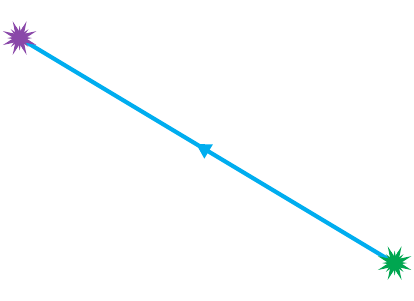}&
\labellist\tiny
\pinlabel \rotatebox{-45}{{\color{my_green} $[2^{[2]},-3]$}} at 35 75
\pinlabel{{\color{my_green} $[3,-2]$}} at 120 55
\pinlabel{{\color{my_purple} $[4,2]$}} at 180 55
\pinlabel{$\frac{2}{7}$} at 7 14
\pinlabel{{\color{my_red} $2$}} at 160 82
\pinlabel{{\color{my_orange} $3$}} at 70 97
\endlabellist
\includegraphics[scale=.5]{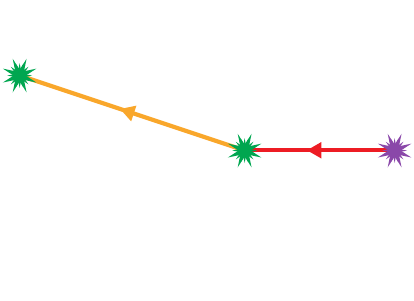}&
\labellist\tiny
\pinlabel \rotatebox{45}{{\color{my_green} $[-2,-3,-2]$}} at 48 106
\pinlabel \rotatebox{45}{{\color{my_purple} $[2,(-2)^{[2]}]$}} at 118 106
\pinlabel \rotatebox{-45}{{\color{my_purple} $[(-2)^{[2]},2]$}} at 117 39
\pinlabel{{\color{my_green} $[3^{[2]}]$}} at 180 18
\pinlabel{{\color{my_red} $2$}} at 52 61
\pinlabel{{\color{my_orange} $3$}} at 147 61
\pinlabel{$\frac{3}{8}$} at 7 14
\endlabellist
\includegraphics[scale=.5]{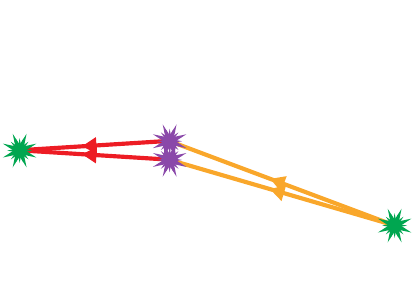}
\\ \hline 
6 & \labellist\tiny
\pinlabel{{\color{my_purple} $[6]$}} at 225 28
\pinlabel{{\color{my_darkblue} $6$}} at 125 85
\pinlabel \rotatebox{-45}{{\color{my_purple} $[(-2)^{[5]}]$}} at 33 119
\pinlabel{$\frac{1}{6}$} at 7 14
\endlabellist
\includegraphics[scale=.5]{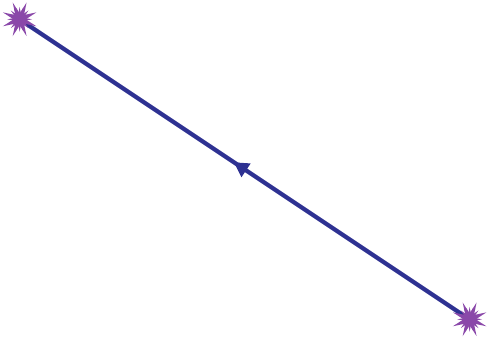}&
\labellist\tiny
\pinlabel{{\color{my_green} $[5,2]$}} at 220 28
\pinlabel{{\color{my_purple} $[4,-2]$}} at 155 28
\pinlabel \rotatebox{0}{{\color{my_green} $[(-2)^{[3]},-3]$}} at 45 133
\pinlabel{{\color{my_green} $4$}} at 93 84
\pinlabel{{\color{my_red} $2$}} at 195 55
\pinlabel{$\frac{2}{9}$} at 7 14
\endlabellist
\includegraphics[scale=.5]{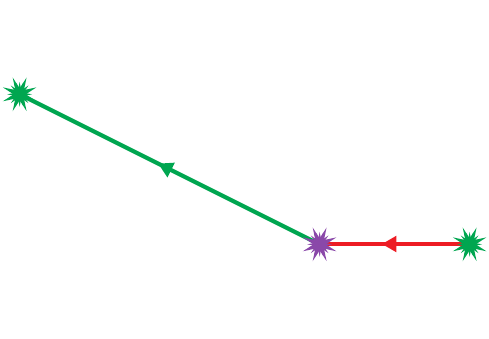}&
\labellist\tiny
\pinlabel{{\color{my_green} $[4,3]$}} at 215 28
\pinlabel\rotatebox{-45}{{\color{my_green} $[3,(-2)^{[2]}]$}} at 145 40
\pinlabel\rotatebox{0}{{\color{my_purple} $[(-2)^{[3]},2]$}} at 120 135
\pinlabel{{\color{my_green} $4$}} at 165 84
\pinlabel{{\color{my_orange} $3$}} at 68 88
\pinlabel{{\color{my_orange} $3$}} at 175 53
\pinlabel{{\color{my_red} $2$}} at 52 125
\pinlabel\rotatebox{-45}{{\color{my_green} $[(-2)^{[2]},-3,-2]$}} at 55 66
\pinlabel{$\frac{3}{11}$} at 12 14
\endlabellist
\includegraphics[scale=.5]{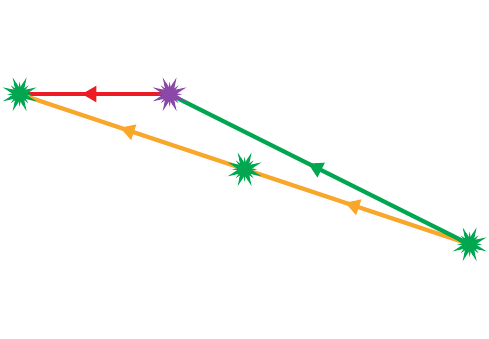}\\[6pt]
 & \labellist\tiny
\pinlabel\rotatebox{0}{{\color{my_purple} $[4,2^{[2]}]$}} at 205 126
\pinlabel{{\color{my_green} $[-3,3]$}} at 110 58
\pinlabel\rotatebox{0}{{\color{my_purple} $[(-2)^{[2]},-4]$}} at 50 126
\pinlabel{{\color{my_orange} $3$}} at 72 96
\pinlabel{{\color{my_orange} $-3$}} at 170 99
\pinlabel{$\frac{3}{10}$} at 12 14
\endlabellist
\includegraphics[scale=.5]{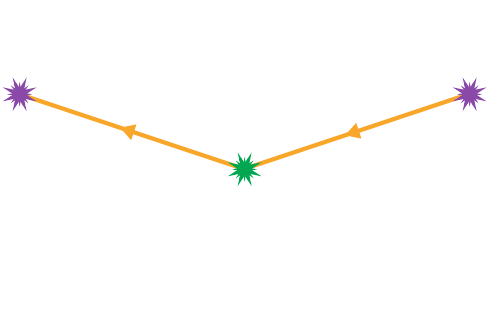}&
\labellist\tiny
\pinlabel\rotatebox{-45}{{\color{my_green} $[3^{[2]},2]$}} at 210 104
\pinlabel\rotatebox{0}{{\color{my_purple} $[(-2)^{[2]},2^{[2]}]$}} at 120 121
\pinlabel\rotatebox{-45}{{\color{my_green} $[3,2,-2]$}} at 180 39
\pinlabel\rotatebox{-45}{{\color{my_green} $[2,-2,-3]$}} at 113 39
\pinlabel\rotatebox{-45}{{\color{my_green} $[-2,(-3)^{[2]}]$}} at 40 39
\pinlabel{{\color{my_orange} $-3$}} at 64 100
\pinlabel{{\color{my_orange} $3$}} at 180 96
\pinlabel{{\color{my_red} $2$}} at 51 63
\pinlabel{{\color{my_red} $2$}} at 125 63
\pinlabel{{\color{my_red} $2$}} at 196 63
\pinlabel{$\frac{5}{13}$} at 12 14
\endlabellist
\includegraphics[scale=.5]{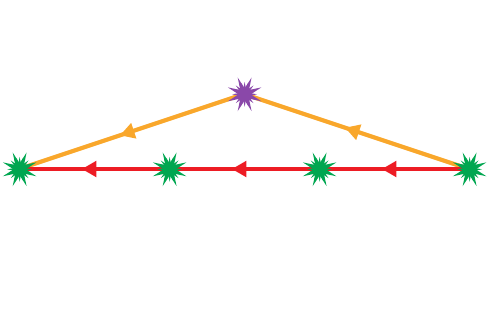}&
\labellist\tiny
\pinlabel \rotatebox{0}{{\color{my_green} $[3,2,3]$}} at 205 60
\pinlabel \rotatebox{45}{{\color{my_purple} $[2,-2,2]$}} at 175 116
\pinlabel\rotatebox{-45}{{\color{my_green} $[-2,-3,2]$}} at 110 40
\pinlabel \rotatebox{45}{{\color{my_green} $[2,-3,-2]$}} at 110 121
\pinlabel \rotatebox{-45}{{\color{my_purple} $[-2,-4,-2]$}} at 45 40
\pinlabel{{\color{my_red} $2$}} at 50 93
\pinlabel{{\color{my_red} $2$}} at 125 93
\pinlabel{{\color{my_red} $2$}} at 197 90
\pinlabel{$\frac{5}{12}$} at 12 14
\endlabellist
\includegraphics[scale=.5]{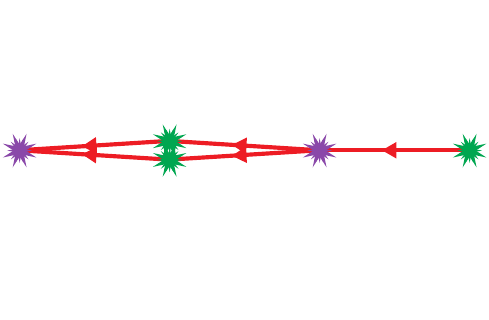}\\ \hline\hline
 \end{tabular}
 \caption{Graphs for 2-bridge knots and links with two through six crossings, with the integer summands $r$ omitted as in Remark \ref{R:integer_part}}
\label{Table:0}
\end{table}

Each vertex in each graph in Tables \ref{Table:0} and \ref{Table:2} corresponds to an essential continued fraction expansion and is labeled accordingly. In each equivalence class of rational number, we choose a representative $\frac{p}{q}$ in the interval $(0,1)$, so the integer summand of each expansion of $\frac{p}{q}$ is nonzero, namely 1, precisely when the first term of the expansion is negative; for compactness, we omit these from our labels---for example, in the first entry in Table \ref{Table:0}, the label $[-2]$ is shorthand for $1+[-2]$.  See Remark \ref{R:integer_part}.

Expansions in which every term is even are marked with purple vertices; the other vertices are green.  Notice that $G_{p/q}$ has two purple vertices, corresponding to the two all-even expansions of $\frac{p}{q}$, whenever $q$ is even.  When $q$ is odd, there is just one such vertex.

 \begin{table}[h!]
\centering
\begin{tabular}{|| c  c  c ||}
 \hline
 \hline
\labellist\tiny
\pinlabel{{\color{my_green} $[7]$}} at 262 25
\pinlabel\rotatebox{-45}{{\color{my_purple} $[(-2)^{[6]}]$}} at 30 155
\pinlabel{{\color{my_purple} $7$}} at 137 85
\pinlabel{$\frac{1}{7}$} at 7 14
\endlabellist
\includegraphics[scale=.5]{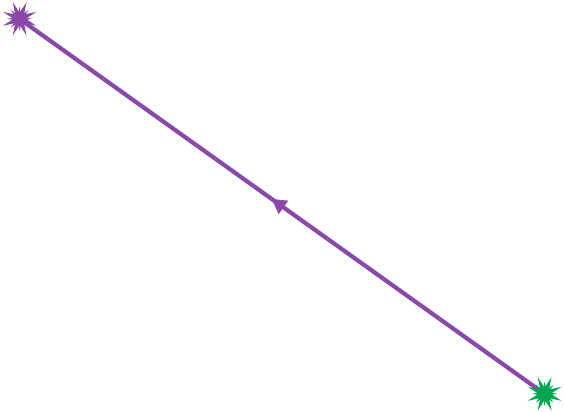}&
\labellist\tiny
\pinlabel{{\color{my_purple} $[6,2]$}} at 250 37
\pinlabel\rotatebox{0}{{\color{my_green} $[5,-2]$}} at 190 37
\pinlabel\rotatebox{-45}{{\color{my_green} $[(-2)^{[4]},-3]$}} at 47 115
\pinlabel{{\color{my_cyan} $5$}} at 108 115
\pinlabel{{\color{my_red} $2$}} at 232 64
\pinlabel{$\frac{2}{11}$} at 12 14
\endlabellist
\includegraphics[scale=.5]{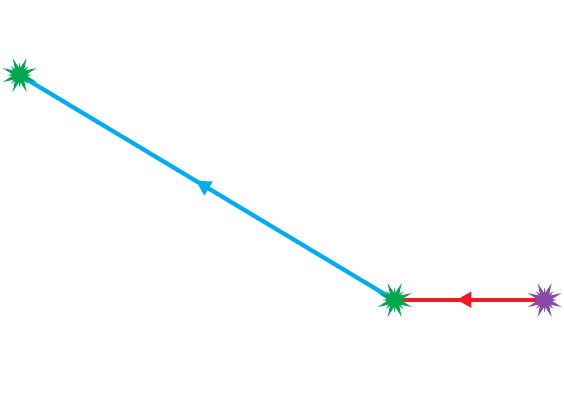}&\\[6pt]
\labellist\tiny
\pinlabel\rotatebox{45}{{\color{my_green} $[5,2^{[2]}]$}} at 245 61
\pinlabel{{\color{my_purple} $[4,-3]$}} at 155 38
\pinlabel\rotatebox{-45}{{\color{my_green} $[(-2)^{[3]},-4]$}} at 42 85
\pinlabel{{\color{my_green} $4$}} at 85 80
\pinlabel{{\color{my_orange} $-3$}} at 208 82
\pinlabel{$\frac{3}{13}$} at 12 14
\endlabellist
\includegraphics[scale=.5]{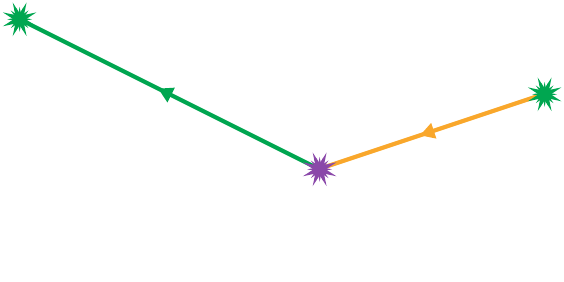}&
\labellist\tiny
\pinlabel\rotatebox{-45}{{\color{my_green} $[(-2)^{[2]},-3,(-2)^{[2]}]$}} at 65 64
\pinlabel\rotatebox{0}{{\color{my_green} $[(-2)^{[3]},3]$}} at 165 103
\pinlabel\rotatebox{-45}{{\color{my_green} $[3,(-2)^{[3]}]$}} at 145 48
\pinlabel{{\color{my_purple} $[4,4]$}} at 235 13
\pinlabel{{\color{my_green} $4$}} at 200 60
\pinlabel{{\color{my_orange} $3$}} at 72 116
\pinlabel{$\frac{4}{15}$} at 12 14
\endlabellist
\includegraphics[scale=.5]{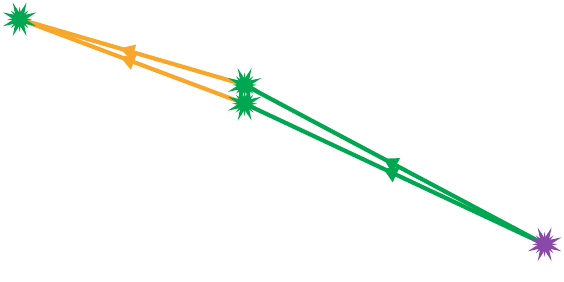}&\\[6pt]
\labellist\tiny
\pinlabel\rotatebox{-45}{{\color{my_green} $[-2,(-3)^{[2]},-2]$}} at 55 47
\pinlabel\rotatebox{-45}{{\color{my_green} $[2,-2,-3,-2]$}} at 110 45
\pinlabel\rotatebox{0}{{\color{my_green} $[-2,-3,-2,2]$}} at 55 114
\pinlabel\rotatebox{-45}{{\color{my_green} $[3,2,-2^{[2]}]$}} at 180 52
\pinlabel\rotatebox{0}{{\color{my_purple} $[2,-2^{[2]},2]$}} at 215 105
\pinlabel\rotatebox{0}{{\color{my_green} $[-2^{[2]},2,3]$}} at 185 125
\pinlabel{{\color{my_green} $[3^{[3]}]$}} at 253 47
\pinlabel{{\color{my_orange} $3$}} at 212 65
\pinlabel{{\color{my_red} $2$}} at 50 88
\pinlabel{{\color{my_red} $2$}} at 123 83
\pinlabel{{\color{my_red} $2$}} at 123 115
\pinlabel{$\frac{8}{21}$} at 12 14
\endlabellist
\includegraphics[scale=.5]{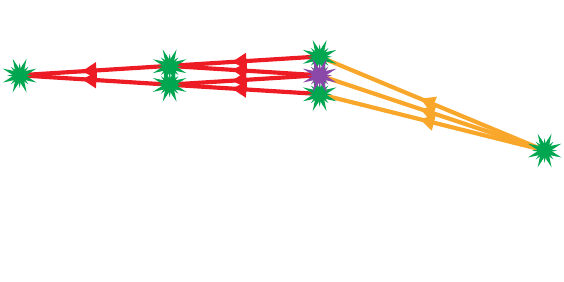}&
\labellist\tiny
\pinlabel \rotatebox{-45}{{\color{my_green} $[-2,-3,-2,-3]$}} at 50 51
\pinlabel{{\color{my_purple} $[(-2)^{[2]},2,-2]$}} at 57 120
\pinlabel \rotatebox{0}{{{\color{my_green} $[(-2)^{[2]},3,2]$}}} at 195 120
\pinlabel \rotatebox{-45}{{{\color{my_green} $[3^{[2]},-2]$}}} at 210 30
\pinlabel\rotatebox{-45}{{\color{my_green} $[2,(-2)^{[2]},-3]$}} at 115 47
\pinlabel \rotatebox{-45}{{{\color{my_green} $[3,4,2]$}}} at 245 88
\pinlabel{{\color{my_orange} $3$}} at 140 68
\pinlabel{{\color{my_orange} $3$}} at 215 89
\pinlabel{{\color{my_red} $2$}} at 51 89
\pinlabel{{\color{my_red} $2$}} at 125 110
\pinlabel{{\color{my_red} $2$}} at 232 55
\pinlabel{$\frac{7}{19}$} at 12 14
\endlabellist
\includegraphics[scale=.5]{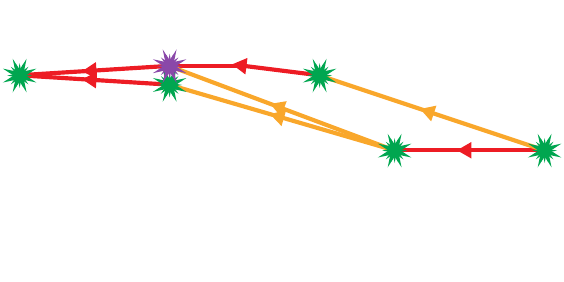}&
\labellist\tiny
\pinlabel\rotatebox{-45}{{\color{my_purple} $[(-2)^{[2]},-4,-2]$}} at 55 58
\pinlabel\rotatebox{0}{{\color{my_green} $[(-2)^{[2]},-3,2]$}} at 140 120
\pinlabel\rotatebox{-45}{{\color{my_green} $[3,-3,-2]$}} at 125 40
\pinlabel \rotatebox{-45}{{{\color{my_green} $[3,-2,2]$}}} at 200 38
\pinlabel \rotatebox{-45}{{{\color{my_green} $[4,2,3]$}}} at 245 100
\pinlabel{{\color{my_orange} $3$}} at 68 79
\pinlabel{{\color{my_orange} $3$}} at 145 95
\pinlabel{{\color{my_red} $2$}} at 52 116
\pinlabel{{\color{my_red} $2$}} at 160 64
\pinlabel{{\color{my_red} $2$}} at 235 64
\pinlabel{$\frac{5}{17}$} at 12 14
\endlabellist
\includegraphics[scale=.5]{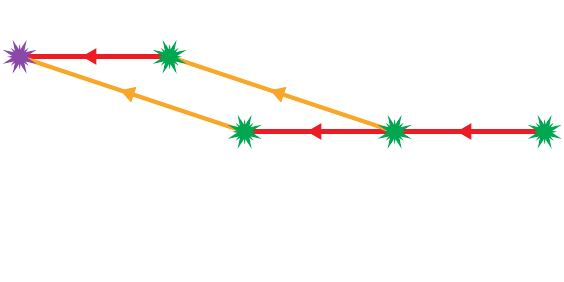}\\[6pt]
\labellist\tiny
\pinlabel\rotatebox{-45}{{\color{my_green} $1+[(-2)^{[3]},-3,-2]$}} at 60 72
\pinlabel \rotatebox{-45}{{{\color{my_purple} $[4,(-2)^{[2]}]$}}} at 160 35
\pinlabel \rotatebox{0}{{{\color{my_purple} $1+[(-2)^{[4]},2]$}}} at 142 142
\pinlabel \rotatebox{-45}{{{\color{my_green} $[5,3]$}}} at 250 50
\pinlabel{$\frac{3}{14}$} at 12 14
\pinlabel{{\color{my_cyan} $5$}} at 180 85
\pinlabel{{\color{my_green} $4$}} at 85 87
\pinlabel{{\color{my_orange} $3$}} at 210 35
\pinlabel{{\color{my_red} $2$}} at 52 144
\endlabellist
\includegraphics[scale=.5]{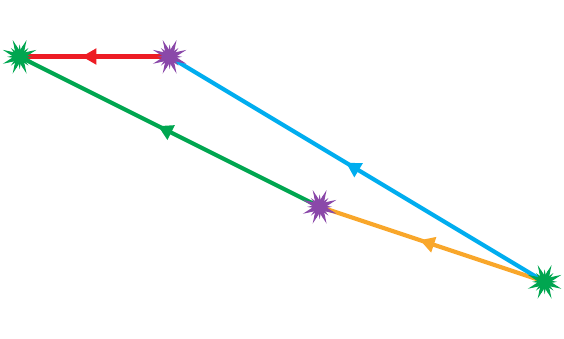}&
\labellist\tiny
\pinlabel \rotatebox{-45}{{\color{my_green} $[(-2)^{[2]},(-3)^{[2]}]$}} at 50 48
\pinlabel \rotatebox{-45}{{{\color{my_green} $[3,-2,-3]$}}} at 152 32
\pinlabel \rotatebox{-45}{{{\color{my_purple} $[4,2,-2]$}}} at 215 30
\pinlabel \rotatebox{0}{{{\color{my_purple} $[(-2)^{[3]},2^{[2]}]$}}} at 172 144
\pinlabel \rotatebox{-45}{{{\color{my_green} $[4,3,2]$}}} at 246 91
\pinlabel{{\color{my_green} $4$}} at 200 103
\pinlabel{{\color{my_orange} $3$}} at 70 88
\pinlabel{{\color{my_orange} $-3$}} at 64 127
\pinlabel{{\color{my_red} $2$}} at 160 55
\pinlabel{{\color{my_red} $2$}} at 235 55
\pinlabel{$\frac{5}{18}$} at 12 14
\endlabellist
\includegraphics[scale=.5]{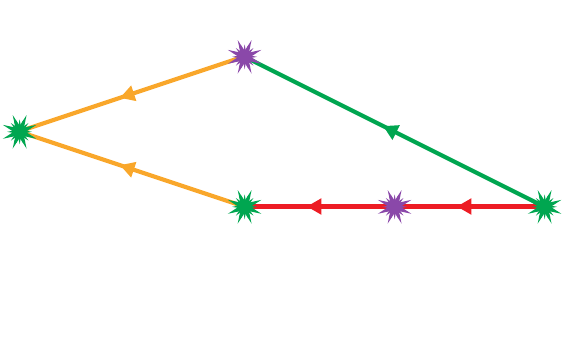}&
\labellist\tiny
\pinlabel \rotatebox{-45}{{\color{my_green} $[-3,(-2)^{[2]},-3]$}} at 52 73
\pinlabel \rotatebox{-45}{{{\color{my_purple} $[2,4,-2]$}}} at 165 126
\pinlabel \rotatebox{-45}{{{\color{my_purple} $[-2,4,2]$}}} at 215 53
\pinlabel \rotatebox{-45}{{\color{my_green} $[-2,3,-2]$}} at 145 53
\pinlabel \rotatebox{-45}{{{\color{my_green} $[2,5,2]$}}} at 245 118
\pinlabel{{\color{my_orange} $3$}} at 70 97
\pinlabel{{\color{my_red} $2$}} at 160 77
\pinlabel{{\color{my_red} $2$}} at 235 80
\pinlabel{$\frac{9}{16}$} at 12 14
\endlabellist
\includegraphics[scale=.5]{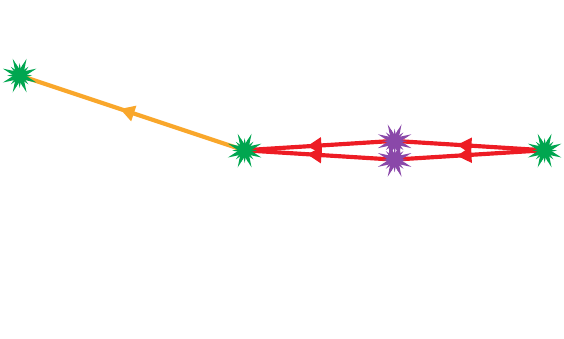}
\\ \hline \hline
 \end{tabular}
\caption{Graphs for 2-bridge knots and links with seven crossings, again with integer summands $r$ mostly omitted from labels as in Remark \ref{R:integer_part}}
\label{Table:2}
\end{table}

Edges in each graph are directed as in Definition \ref{D:G_p/q}. Each graph $G_{p/q}$ is drawn in a way that reflects the geography of the knot or link $K_{p/q}$.  Thus, each edge is directed from right to left, and shorter expansions appear below longer ones.   Vertices that appear to overlap map to the same geography point.  

More precisely, each edge is colored and sloped according to the type of move it signifies. Each red edge has zero slope (except in some cases of overlapping vertices) and indicates a $2$-move.  Each orange edge that points upward (and leftward) corresponds to a 3-move, while each orange edge that points downward indicates an inverse $-3$-move. Each green edge describes a 4-move, each light blue edge a 5-move, each dark blue edge a 6-move, and the purple edge a 7-move.

\begin{remark}
The crossing number of a knot or link $K_{p/q}$ manifests in the graph as the edge-weighted distance between $\mathbf{b}_+$ and $\mathbf{b}_-$.  That is, for any directed path $P$ in $G_{p/q}$ from $\mathbf{b}_+$ to $\mathbf{b}_-$, if $n_1,\dots,n_r$ are the labels on the edges of $P$, then $|n_1|+\cdots+|n_r|$ equals the crossing number of $K_{p/q}$. One way for an expert to see this is to use Theorem 2.1 of \cite{CohKinLowShaVan}, which states that, given any reduced alternating diagram $D$ of this link, every continued fraction $\mathbf{b}_i$ in the sequence $\mathbf{b}_+=\mathbf{b}_0\to\mathbf{b}_1\to\cdots\to\mathbf{b}_r=\mathbf{b}_-$ corresponds to a state surface from $D$; in particular $\mathbf{b}_+$ and $\mathbf{b}_-$ respectively correspond to the all-$B$ and all-$A$ checkerboard surfaces from $D$. Each move $\mathbf{b}_{i-1}\to\mathbf{b}_i$ corresponds to switching $|n_i|$ $B$-smoothings to $A$-smoothings, so the sum $|n_1|+\cdots+|n_r|$ equals the total number of crossings in $D$, which equals the crossing number of $K_{p/q}$ \cite{Kauffman1987,Murasugi1987,Thistlethwaite1987}.
\end{remark}

The next two remarks explain how the preceding tables describe $G_{p'/q'}$ for any rational number which $\frac{p'}{q'}$ is equivalent to one of the numbers $\frac{p}{q}$ in the tables, in the sense that $q'=q$ and either $p'\equiv \pm p\pmod{q}$ or $pp'\equiv \pm 1\pmod{q}$.

\begin{figure}
    \centering
    \labellist\tiny
\pinlabel\rotatebox{-45}{{\color{my_green} $1+[-2,-3,(-2)^{[3]}]$}} at 55 72
\pinlabel \rotatebox{-45}{{{\color{my_purple} $1+[(-2)^{[2]},4]$}}} at 160 33
\pinlabel \rotatebox{0}{{{\color{my_purple} $[2,(-2)^{[4]}]$}}} at 130 142
\pinlabel \rotatebox{-45}{{{\color{my_green} $[3,5]$}}} at 250 50
\pinlabel{$\frac{5}{14}$} at 12 14
\pinlabel{{\color{my_cyan} $5$}} at 180 85
\pinlabel{{\color{my_green} $4$}} at 85 87
\pinlabel{{\color{my_orange} $3$}} at 210 35
\pinlabel{{\color{my_red} $2$}} at 52 144
\endlabellist
\includegraphics[scale=.5]{Figures/Table/53}
\hfill
\labellist\tiny
\pinlabel\rotatebox{0}{{\color{my_green} $1+[-3,-5]$}} at 78 20
\pinlabel \rotatebox{-45}{{{\color{my_purple} $[2^{[2]},-4]$}}} at 163 37
\pinlabel \rotatebox{0}{{{\color{my_purple} $1+[-2,2^{[4]}]$}}} at 155 142
\pinlabel \rotatebox{-45}{{{\color{my_green} $[2,3,2^{[3]}]$}}} at 252 167
\pinlabel{$\frac{9}{14}$} at 12 14
\pinlabel{{\color{my_cyan} $-5$}} at 110 92
\pinlabel{{\color{my_green} $-4$}} at 205 89
\pinlabel{{\color{my_orange} $-3$}} at 93 42
\pinlabel{{\color{my_red} $2$}} at 230 143
\endlabellist
\includegraphics[scale=.5]{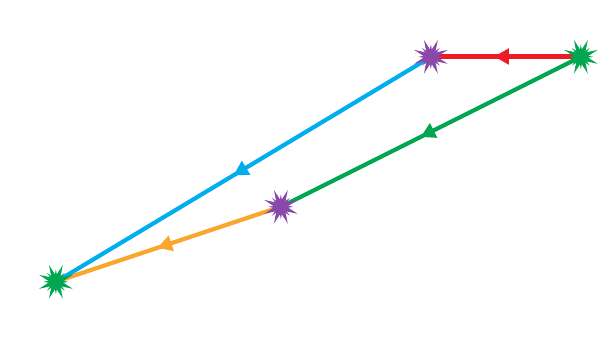}
\hfill
\labellist\tiny
\pinlabel\rotatebox{0}{{\color{my_green} $1+[-5,-3]$}} at 78 20
\pinlabel \rotatebox{-45}{{{\color{my_purple} $1+[-4,2^{[2]}]$}}} at 174 25
\pinlabel \rotatebox{0}{{{\color{my_purple} $[2^{[4]},-2]$}}} at 168 142
\pinlabel \rotatebox{-45}{{{\color{my_green} $[2^{[3]},3,2]$}}} at 252 167
\pinlabel{$\frac{11}{14}$} at 12 14
\pinlabel{{\color{my_cyan} $-5$}} at 110 92
\pinlabel{{\color{my_green} $-4$}} at 205 89
\pinlabel{{\color{my_orange} $-3$}} at 93 42
\pinlabel{{\color{my_red} $2$}} at 230 143
\endlabellist
\includegraphics[scale=.5]{Figures/Table/_5_3}
\caption{The graphs $G_{5/14}$, $G_{9/14}$, and $G_{11/14}$ obtained from $G_{3/14}$ in Table~\ref{Table:2} as in Remarks \ref{R:Negate} and \ref{R:Invert}}
\label{Fi:Other_p_q}
\end{figure}

\begin{remark}\label{R:Negate}
  One can construct the graph $G_{-p/q}$, or equivalently $G_{(q-p)/q}$, from $G_{p/q}$ by reflecting the graph and its labels across a vertical axis, reversing the arrows on the edges, negating each edge label, and negating each term in each vertex label.   For example, the graph $G_{11/14}$ in Figure \ref{Fi:Other_p_q} can be obtained in this way from the graph $G_{3/14}$ in Table \ref{Table:2}.
  Some of the graphs corresponding to knots or links with even crossing number have reflection symmetry across a vertical axis; the associated knots and links are amphichiral.  The knots and links associated to the other graphs are chiral. In these cases, the links $K_{p/q}$ and $K_{(q-p)/q}$ are not isotopic.
\end{remark}

Remark \ref{R:Negate} is particularly noteworthy because, while no edge in Table \ref{Table:0}  
nor \ref{Table:2} describes an inverse move of type $-4$, $-5$, $-6$, nor $-7$,  these moves will appear in several of the mirrored graphs, including $G_{3/4}$, $G_{4/5}$, $G_{5/6}$, $G_{6/7}$.

\begin{remark}\label{R:Invert}
    One basic number theoretic property of continued fractions is that, if $r+[b_1,\dots,b_k]=\frac{p}{q}$ and $r'+[b_k,\dots,b_1]=\frac{p'}{q'}$, then $q=q'$ and $pp'\equiv 1\pmod{q}$. See for example Hatcher \cite[Section 2.1]{Hatcher}. 
    Given such $p,q,p'$, the graphs $G_{p/q}$ and $G_{p'/q}$ are identical, except in their vertex labels. One obtains each vertex label in $G_{p'/q}$ by reversing the order of the terms in the corresponding vertex label in $G_{p/q}$ and then adjusting the integer summand as specified by Remark \ref{R:integer_part}. For example, $G_{5/14}$ in Figure \ref{Fi:Other_p_q} is obtained in this way from $G_{3/14}$ in Table \ref{Table:2}; likewise, either of $G_{9/14}$ or $G_{11/14}$ in Figure \ref{Fi:Other_p_q} can be obtained in this way from the other.
\end{remark}

\subsection{Graphs for infinite families of rational numbers}
\label{subsec:families}

For many infinite families of rational numbers $\frac{p}{q}$, it is straightforward to characterize the edge-weighted directed graph $G_{p/q}$ precisely. This is easiest when $\frac{p}{q}$ has a particularly short continued fraction expansion, so we do this next for all rational numbers with expansions of length at most two.

\begin{figure}[h!]
    \centering
    \;\hfill
    \labellist\tiny
    \pinlabel{$\frac{1}{2}$} at 5 8
    \pinlabel{$2$} at 68 12
    \pinlabel{$[2]$} at 100 28
    \pinlabel{$1+[-2]$} at 30 28
    \endlabellist
    \includegraphics[scale=.8]{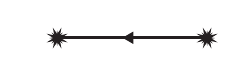}\hfill
    \labellist\tiny
    \pinlabel{$\frac{1}{q}$; $q\geq 3$} at 19 9
    \pinlabel{$q$} at 66 21
    \pinlabel{$[q]$} at 100 30
    \pinlabel{$1+[(-2)^{[q-1]}]$} at 35 47
    \endlabellist
    \includegraphics[scale=.8]{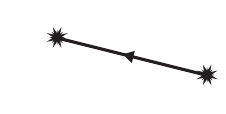}\hfill  
    \labellist\tiny
    \pinlabel{$-q$} at 92 22
    \pinlabel{$\frac{1-q}{q}$; $q\geq 3$} at 15 9
    \pinlabel{$[2^{[q-1]}]$} at 108 45
    \pinlabel{$1+[-q]$} at 39 28
    \endlabellist
    \includegraphics[scale=.8]{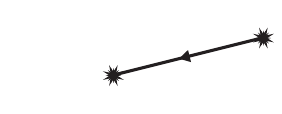}\hfill\;
    \caption{The graphs for the family of rational numbers from Example \ref{Ex:q}}
    \label{Fi:Families1}
\end{figure}

\begin{example}\label{Ex:q}
Consider a rational number of the form $\pm\frac{1}{q}$ for some integer $q\geq 2$.  These correspond to the class of 2-bridge torus knots and links $K_{\pm 1/q}$.  
    The number $\frac{1}{q}$ has precisely two essential continued fraction expansions: $[q]$ and $1+[(-2)^{q-1}]$. Similarly, $\frac{-1}{q}$ is equivalent to $\frac{1-q}{q}$, which has precisely two essential continued fraction expansions: $1+[-q]$ and $[2^{[q-1]}]$. Thus, every graph $G_{\pm 1/q}$ has two vertices and a single edge, which is weighted by $\pm q$. See Figure \ref{Fi:Families1}.
\end{example}

\begin{figure}[h!]
    \centering
    \;\hfill
    \labellist\tiny
    \pinlabel{$\frac{2}{5}$} at 5 9
    \pinlabel{$2$} at 76 21
    \pinlabel{$2$} at 149 21
    \pinlabel{$[3,2]$} at 178 37
    \pinlabel{$[2,-2]$} at 110 37
    \pinlabel{$1+[-2,-3]$} at 30 37
    \endlabellist    \includegraphics[scale=.8]{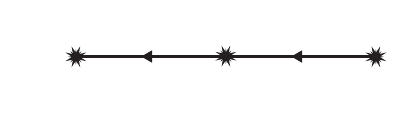}
    \hfill
    \labellist\tiny
    \pinlabel{$\frac{2}{2q+1}$; $q\geq 3$} at 27 9
    \pinlabel{$q$} at 75 21
    \pinlabel{$2$} at 149 13
    \pinlabel{$[q+1,2]$} at 168 28
    \pinlabel{$[q,-2]$} at 113 28
    \pinlabel{$1+[(-2)^{[q-1]},-3]$} at 50 45
    \endlabellist    \includegraphics[scale=.8]{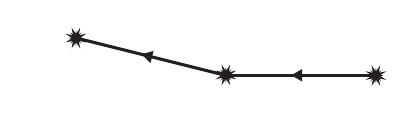}
    \hfill\;
    \caption{The graphs for the family of rational numbers from Example \ref{Ex:k2}}
    \label{Fi:Families2}
\end{figure}

\begin{example}\label{Ex:k2}
Consider a rational number of the form $\frac{2}{2q+1}$ with $q\geq 2$.  These correspond to the class of twist knots $K_{2/(2q+1)}$.  
The positive continued fraction expansion of  $\frac{2}{2q+1}$ is given by $[q+1,2]$.  This expansion admits just one essential $n$-move, namely the $2$-move $[q+1,2]\to[q,-2]$. The resulting expansion $[q,-2]$ admits just two essential moves: the $-2$-move back to $[q+1,2]$ and the $q$-move $[q,-2]\to1+[(-2)^{[q-1]},-3]$, so Theorem \ref{T:essentialconnected} implies that for every integer $q\geq 2$, the graph $G_{2/(2q+1)}$ is a path graph of length two, with one edge weighted by 2 and the other weighted by $q$. See Figure \ref{Fi:Families2}.
\end{example}

\begin{figure}[h!]
    \centering
    \;\hfill
    \labellist\tiny
    \pinlabel{$[3,3]$} at 150 26
    \pinlabel \rotatebox{-45}{{$1+[-2^{[2]},2]$}} at 100 31
    \pinlabel \rotatebox{0}{{$[2,(-2)^{[2]}]$}} at 105 80
    \pinlabel \rotatebox{-45}{{$1+[-2,-3,-2]$}} at 38 32
    \pinlabel {$\frac{3}{8}$} at 5 9
    \pinlabel{$2$} at 50 72
    \pinlabel{$3$} at 123 55
    \endlabellist    \includegraphics[scale=.8]{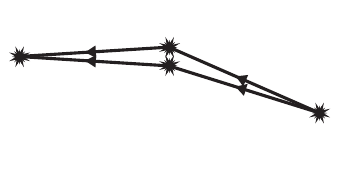}
    \hfill
    \labellist\tiny
    \pinlabel{$[q,3]$} at 150 26
    \pinlabel \rotatebox{-30}{{$[q-1,(-2)^{[2]}]$}} at 90 30
    \pinlabel \rotatebox{0}{{$1+[(-2)^{[q-1]},2]$}} at 127 83
    \pinlabel \rotatebox{-45}{{$1+[-2,-3,-2]$}} at 34 48
    \pinlabel {$\frac{3}{3q-1}$; $q\geq 4$} at 28 9
    \pinlabel{$2$} at 50 76
    \pinlabel{$q$} at 124 61
    \pinlabel{$3$} at 121 34
    \pinlabel{$q-1$} at 61 66
    \endlabellist    \includegraphics[scale=.8]{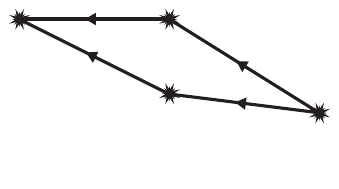}\hfill
    \labellist\tiny
    \pinlabel{$[q,\ell]$} at 150 26
    \pinlabel \rotatebox{-30}{{$[q-1,(-2)^{[\ell-1]}]$}} at 100 25
    \pinlabel \rotatebox{0}{{$1+[(-2)^{[q-1]},\ell-1]$}} at 113 83
    \pinlabel \rotatebox{0}{{$\mathbf{b}_-$}} at 8 72
    \pinlabel {$\frac{\ell}{\ell q-1}$; $q,\ell\geq 4$} at 31 9
    \pinlabel{$\ell-1$} at 40 83
    \pinlabel{$q$} at 124 57
    \pinlabel{$\ell$} at 121 34
    \pinlabel{$q-1$} at 39 54
    \endlabellist    \includegraphics[scale=.8]{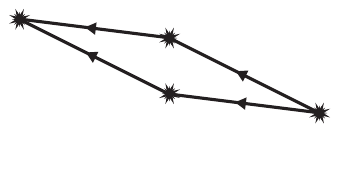}\hfill
    \caption{The graphs  for the families of rational number from Example \ref{Ex:kl}. On the far right, $\mathbf{b}_-=1+[(-2)^{[q-2]},-3,(-2)^{[\ell-2]}]$. }
    \label{Fi:Families3}
\end{figure}

\begin{example}\label{Ex:kl}
    Given a rational number of the form $\frac{3}{3q-1}$ where $q\geq 3$, its positive expansion is $[q,3]$.  By similar reasoning as in Example \ref{Ex:k2}, Theorem \ref{T:essentialconnected} implies that $G_{3/(3q-1)}$ is a cycle graph whose four edges are weighted $q$, $3$, $q-1$, and $2$, as shown left and center in Figure \ref{Fi:Families3}. More generally, given any rational number of the form $\frac{\ell}{q\ell-1}=
    [q,\ell]$ for some $q,\ell\geq 3$, the graph $G_{\ell/(q\ell-1)}$ is a cycle graph whose four edges are weighted $q$, $\ell$, $q-1$, and $\ell-1$. Figure \ref{Fi:Families3} depicts all three cases, depending on how many of $q$ and $\ell$ equal 3.
\end{example}

\begin{figure}[h!]
    \centering
\;\hfill
    \labellist\tiny
    \pinlabel{$[q+1,2^{[2]}]$} at 155 54
    \pinlabel \rotatebox{0}{{$[q,-3]$}} at 100 18
    \pinlabel \rotatebox{0}{{$1+[(-2)^{[q-1]},-4]$}} at 45 54
    \pinlabel {$\frac{3}{3q+1}$; $q\geq 3$} at 28 9
    \pinlabel{$-3$} at 140 30
    \pinlabel{$q$} at 56 31
    \endlabellist   \includegraphics[scale=.8]{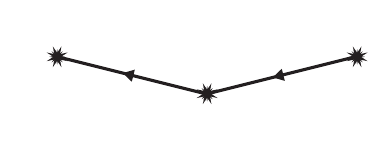}
    \hfill
    \labellist\tiny
    \pinlabel{$[q+1,2^{[\ell-1]}]$} at 148 54
    \pinlabel \rotatebox{0}{{$[q,-\ell]$}} at 100 18
    \pinlabel \rotatebox{0}{{$1+[(-2)^{[q-1]},-\ell-1]$}} at 53 54
    \pinlabel {$\frac{\ell}{\ell q+1}$; $q,\ell\geq 3$} at 32 9
    \pinlabel{$-\ell$} at 140 30
    \pinlabel{$q$} at 56 31
    \endlabellist    \includegraphics[scale=.8]{Figures/k_l}
    \hfill\;

    \caption{The graphs for the families of rational numbers  from Example \ref{Ex:q-l}}
    \label{Fi:Families4}
\end{figure}

\begin{example}\label{Ex:q-l}
    Every rational number of the form $\frac{3}{3q+1}$ where $q\geq 3$ has an expansion of the form $[q,-3]$. Its positive expansion is thus $[q+1,2,2]$, and Theorem \ref{T:essentialconnected} implies that $G_{3/(3q+1)}$ is a path graph whose two edges are weighted $-3$, and $q$, as shown on the left in Figure \ref{Fi:Families4}.  More generally, given any rational number of the form $\frac{\ell}{q\ell+1}=
[q,-\ell]$ for some $q,\ell\geq 3$, the graph $G_{\ell/(q\ell+1)}$ is a path graph whose two edges are weighted $q$ and $-\ell$, as shown in Figure \ref{Fi:Families4}. Note in particular that when $\ell=q$ the positive and negative expansions $\mathbf{b}_+$ and $\mathbf{b}_-$ of the fraction $\frac{q}{q^2+1}$ have the same length. More on this shortly.
\end{example}

We conclude by highlighting a noteworthy property of a certain subfamily of the rational numbers from Example \ref{Ex:q-l}. 
Given $\frac{p}{q}\in\mathbb{Q}\setminus\mathbb{Z}$, let $I^\text{min}_{p/q}$ denote the induced subgraph of $G_{p/q}$ associated with the minimal length expansions of $\frac{p}{q}$. Theorem~4 in Hirasawa-Teragaito \cite{Hirasawa_Teragaito:2006} implies that $I^\text{min}_{p/q}$ is connected, as all minimal length expansions of $\frac{p}{q}$ are related by essential $\pm 2$-moves. The following example, however, demonstrates that the analogous result for the induced subgraph  $I^\text{max}_{p/q}$ of $G_{p/q}$ associated with the {\it maximal} length expansions of $\frac{p}{q}$ is false. 

\begin{example}
    \label{Ex:max}
    As shown in Table \ref{Table:0}, the rational number $\frac{p}{q}=\frac{3}{10}$ has exactly two essential expansions of maximal length, namely $[4,2^{[2]}]$ and $1+[(-2)^{[2]},-4]$, and these expansions are not directly related by any $\pm n$-move. The induced subgraph $I^\text{max}_{3/10}$ is therefore disconnected.
    
    More generally, given any integer $q\geq 3$, there are precisely two maximal length expansions of ${q/(q^2+1)}=[q,-q]$, namely $1+[(-2)^{[q-1]},-(q+1)]$ and $[q+1,2^{[q-1]}]$, and these expansions are not directly related by any $\pm n$ move.  Hence, for each integer $q\geq 3$, the graph $I^\text{max}_{q/(q^2+1)}$ is disconnected.
\end{example}

To recapitulate, the preceding example establishes the following result. 

\begin{proposition}
    \label{prop:maximal}
The induced subgraph  $I^\text{max}_{p/q}$ of $G_{p/q}$ associated with the {\it maximal} length expansions of $\frac{p}{q}$ need not be connected.  In particular, $I^\text{max}_{p/q}$ is disconnected whenever $\frac{p}{q}=\frac{q}{q^2+1}=[q,-q]$ for any $q\geq 3$.
\end{proposition}

\appendix
\section{Full list of moves}
\label{allmoves}

All moves on continued fraction expansions discussed in this paper are described here, including explicit formulas for all edge cases. 
Throughout this Appendix, $n$ denotes an integer greater than or equal to $2$ and $\varepsilon=\pm1$.

\begin{itemize}[leftmargin=*, widest=\textbf{$-n$-move and inverse: }]
\item[\textbf{0-move: }] Removing a term $b_i=0$
        \begin{align*}
            r+[\ldots,a,0,b,\ldots] &\rightarrow r+[\ldots,a+b,\ldots],\\
            r+[\ldots,a,0] &\rightarrow r+[\ldots,a],\\
      r+[0,a,b,\ldots] &\rightarrow (r-a)+[b,\ldots].
\end{align*} 

\item[\textbf{$\boldsymbol{\pm 1}$-move: }]
Removing a term $b_i=\upvarepsilon=\pm1$
        \begin{align*}
            r+[\ldots,a,\upvarepsilon,b,\ldots] &\rightarrow r+[\ldots,a-\upvarepsilon,b-\upvarepsilon,\ldots],\\
            r+ [\ldots,a,\upvarepsilon] &\rightarrow r+ [\ldots,a-\upvarepsilon],\\
            r+[\upvarepsilon,a,\ldots]&\rightarrow (r+\upvarepsilon)+[a-\upvarepsilon,\ldots].
        \end{align*}
        
\item[\textbf{232-move: }]
Replacing a subsequence of the form $2\upvarepsilon, (3\upvarepsilon)^{[n-2]}, 2\upvarepsilon$ with $(-3\upvarepsilon)^{[n-1]}$
        \begin{align*}
            r+[\ldots,a,2\upvarepsilon, (3\upvarepsilon)^{[n-2]}, 2\upvarepsilon,b,\ldots] &\rightarrow r+[\ldots,a-\upvarepsilon,(-3\upvarepsilon)^{[n-1]}, b-\upvarepsilon,\ldots],\\
            r+[\ldots,a,2\upvarepsilon, (3\upvarepsilon)^{[n-2]}, 2\upvarepsilon] &\rightarrow r+[\ldots,a-\upvarepsilon,(-3\upvarepsilon)^{[n-1]}],\\
            r+[2\upvarepsilon, (3\upvarepsilon)^{[n-2]}, 2\upvarepsilon,a,\ldots]&\rightarrow (r+\upvarepsilon)+[(-3\upvarepsilon)^{[n-1]}, a-\upvarepsilon,\ldots],\\
            r+[2\upvarepsilon, (3\upvarepsilon)^{[n-2]}, 2\upvarepsilon]&\rightarrow (r+\upvarepsilon)+[(-3\upvarepsilon)^{[n-1]}].
        \end{align*}

\item[\textbf{$\boldsymbol{\pm2}$-move: }] Replacing a term $b_i=2 \varepsilon $ with $b_i=-2 \varepsilon $ where $\varepsilon = \pm 1$
        \begin{align*}
            r+[\ldots,a,2\varepsilon ,b,\ldots] &\rightarrow r+[\ldots,a-\varepsilon,-2\varepsilon, b-\varepsilon,\ldots],\\
            r+[\ldots,a,2\varepsilon] &\rightarrow r+[\ldots,a-\varepsilon,-2\varepsilon],\\
            r+[2\varepsilon, a,\ldots] &\rightarrow (r+\varepsilon)+[-2\varepsilon, a-\varepsilon,\ldots],\\
            r+[2\varepsilon]&\rightarrow (r+\varepsilon)+[-2\varepsilon].
        \end{align*}

\item[\textbf{$\pm n$-move and inverse: }] Replacing a term $b_i=n\varepsilon$ with a subsequence $(-2\varepsilon)^{[n-1]}$ and vice versa for $n\geq2$, $\varepsilon=\pm1$
        \begin{align*}
            r+[b_1,\dots,b_{i-1},n\varepsilon,b_{i+1},\dots,b_k]&\leftrightarrow r+[b_1,\dots,b_{i-1}-\varepsilon,(-2\varepsilon)^{[n-1]},b_{i+1}-\varepsilon,\dots,b_k],\\
            r+[b_1,\dots,b_{k-1},n\varepsilon]&\leftrightarrow r+[b_1,\dots,b_{k-1}-\varepsilon,(-2\varepsilon)^{[n-1]}],\\
            r+[n\varepsilon,b_{2},\dots,b_k]&\leftrightarrow (r+\varepsilon)+[(-2\varepsilon)^{[n-1]},b_{2}-\varepsilon,\dots,b_k],\\
            r+[n\varepsilon]&\leftrightarrow (r+\varepsilon)+[(-2\varepsilon)^{[n-1]}].
        \end{align*}
\end{itemize}

\bibliographystyle{plain}
%\bibliography{References}

\end{document}